\documentclass[10pt]{amsart}

\usepackage{latexsym,exscale,enumerate,amsfonts,amssymb,xparse, mathtools}
\usepackage[normalem]{ulem}
\usepackage{amsmath,amsthm,amsfonts,amssymb,amscd, stmaryrd,textcomp,bbm,euscript}
\usepackage[bookmarks=false]{hyperref}
\usepackage{bbold}
\usepackage[usenames,dvipsnames]{xcolor}
\usepackage{enumitem}
\usepackage{verbatim}
\usepackage{ytableau}

\usepackage{amsmath,semantex}

\NewSymbolClass\MyBinaryOperator[
    define keys={
        {Lder}{command=\overset{L}},
        {Rder}{upper=R},
    },
]
\NewObject\MyBinaryOperator\tensor{\otimes}[
    define keys={
        {der}{Lder},
    },
]
\NewObject\MyBinaryOperator\fibre{\times}[
    define keys={
        {der}{Rder},
    },
]
\allowdisplaybreaks
\usepackage{array}
\usepackage{scalerel}

\definecolor{ourblue}{RGB}{109, 156, 179}

\usepackage{xcolor}
\definecolor{myorange}{rgb}{1,0.647,0}

\definecolor{mypurple}{cmyk}{.6,.9,0, .11}

\numberwithin{equation}{section}

\theoremstyle{definition}
\newtheorem{thm}{Theorem}[section]
\newtheorem{cor}[thm]{Corollary}

\newtheorem{conj}[thm]{Conjecture}
\newtheorem{lem}[thm]{Lemma}
\newtheorem{rem}[thm]{Remark}

\newtheorem{prop}[thm]{Proposition}
\newtheorem{defn}[thm]{Definition}
\newtheorem{ex}[thm]{Example}

\newtheorem{notn}[thm]{Notation}

\usepackage{tikz}
\usetikzlibrary{cd}
\usetikzlibrary{snakes}
\usetikzlibrary{decorations.markings}
\usetikzlibrary{decorations.pathreplacing}
\usetikzlibrary{arrows,shapes,positioning}
\usetikzlibrary{decorations.pathmorphing}

\tikzset{anchorbase/.style={baseline={([yshift=-0.5ex]current bounding box.center)}},
tinynodes/.style={font=\tiny,text height=0.75ex,text depth=0.15ex},
smallnodes/.style={font=\scriptsize,text height=0.75ex,text depth=0.15ex},
>={Latex[length=1mm, width=1.5mm]},
overcross/.style={line width=4pt,color=white},
1label/.style={very thick, black},
2label/.style={very thick, ourblue},
m1label/.style={very thick, double, black},
m2label/.style={very thick, double, ourblue},
glabel/.style={very thick, gray},
mglabel/.style={very thick, double, gray},
Klabel/.style={line width=5pt, mypurple},
webs/.style={line width=.9,color=black},
cover/.style={line width=5pt,color=white}
}

\tikzstyle directed=[postaction={decorate,decoration={markings,
    mark=at position #1 with {\arrow{>}}}}]
\tikzstyle rdirected=[postaction={decorate,decoration={markings,
    mark=at position #1 with {\arrow{<}}}}]

\newcommand{\Hom}{{\rm Hom}}

\renewcommand{\to}{\rightarrow}

\newcommand{\Cone}{{\rm Cone}}

\let\hat=\widehat
\let\tilde=\widetilde
\def\C{{\mathbb C}}

\def\Z{{\mathbb Z}}

\def\X{{\mathbb{X}}}
\def\Y{{\mathbb{Y}}}
\def\M{{\mathbb{M}}}
\def\F{{\mathbb{F}}}
\def\A{{\mathbb{A}}}
\def\B{{\mathbb{B}}}
\def\deg{{\mathrm{deg}}}
\def\p{{\tilde{p}}}

\renewcommand{\L}{\mathcal{L}}

\newcommand{\SG}{\mathfrak{S}}

\newcommand{\Sym}{\mathrm{Sym}}

\newcommand{\CS}{\EuScript C}

\newcommand{\RS}{\mathcal{R}}
\newcommand{\IS}{\mathcal{I}}
\newcommand{\JS}{\mathcal{J}}
\newcommand{\VS}{\mathcal{V}}
\newcommand{\Scal}{\mathcal{S}}
\newcommand{\Bim}{\mathrm{Bim}}
\newcommand{\Mod}{\mathrm{Mod}}

\newcommand{\SSBim}{\mathrm{SSBim}}

\renewcommand{\aa}{\mathbf{a}}
\newcommand{\bb}{\mathbf{b}}
\newcommand{\cc}{\mathbf{c}}

\definecolor{revisions}{RGB}{0,0,0}

\begin{document}
%

\title{Doubly-Graded Exponential Growth of Colored Torus Knot Homology}
\author{Luke Conners}


\address{Institute of Mathematics, Universität Zürich, Wintherthurerstrasse 190, CH-0857 Zürich}
\email{luke.conners@math.uzh.ch}

\begin{abstract}
We give an invariant construction of reduced HOMFLY homology for arbitrary links reduced at components of arbitrary color and prove some structural properties relating this invariant to unreduced HOMFLY homology. Combined with previous results, this gives a recursive formula for the reduced HOMFLY homology of colored positive torus knots and some colored positive torus links. We prove that after forgetting the quantum grading, the resulting doubly-graded invariant of positive torus knots grows exactly exponentially in the color, resolving a 2013 conjecture of Gorsky--Gukov--Stošić. Finally, we verify the ``color-shifting" conjecture of Gukov--Nawata--Saberi--Stošić--Sułkowski in many examples.
\end{abstract}

\maketitle

\setcounter{tocdepth}{1}
\tableofcontents

\section{Introduction}

\subsection{Context and Statement of Results} \label{subsec: intro}

In \cite{KhR08}, Khovanov and Rozansky defined a triply-graded homological invariant of links called \textit{HOMFLY homology}, categorifying the famous HOMFLY polynomial. Somewhat surprisingly, this invariant was conjectured to exist by Dunfield--Gukov--Rasmussen in \cite{DGR06} before its rigorous construction. From the perspective of those authors, the existence of HOMFLY homology was strongly suggested by the rigid structure underlying the various existing $\mathfrak{sl}(N)$ link homologies, also due to Khovanov--Rozansky \cite{KhR08I}. The DGR conjectures went far beyond merely the existence of HOMFLY homology, predicting a wealth of structural properties which have driven a thriving research program for over 20 years. Most notably, we mention the spectral sequences to $\mathfrak{gl}(N)$ homology \cite{Ras15}, to a ``symmetric" $\mathfrak{gl}(-N)$ homology \cite{RW19, RW20, QRS18}, and (very recently) to knot Floer homology \cite{BPRW22}; and the ``mirror symmetry" conjecture, resolved by Oblomkov--Rozansky and Gorsky--Hogancamp--Mellit using deep and surprising connections with the algebraic geometry of Hilbert schemes \cite{GHM21, GH22, OR18, OR19b, OR20, GNR21}. Determining when these spectral sequences collapse and unraveling the depth of this geometric construction remain active areas of research even 20 years later.

HOMFLY homology admits several generalizations to invariants of \textit{colored} links, in which each link component is labeled by a positive integer. Following their rigorous construction by Mackaay--Stošić--Vaz and Webster--Williamson \cite{MSV11, WW17}, many structural properties of these invariants were also conjectured by Gukov--Stošić in \cite{GS12}. To discover their conjectures, those authors rely on a physical interpretation of link homology, through the enumerative geometry of BPS states arising from string theory, established by Gukov--Schwarz--Vafa in \cite{GSV05}. This physical connection has proved to be an especially rich perspective, generating a wealth of new conjectures surrounding colored HOMFLY homology by Gorsky--Gukov--Stošić and Gukov--Nawata--Saberi--Stošić--Sułkowski in \cite{GGS18, GNSSS16}. For convenience, we refer to these as the \textit{GGS} and \textit{GNSSS} conjectures.

This paper is concerned with the investigation of the GGS and GNSSS conjectures. To formulate their conjectures, GGS and GNSSS rely on a conjectural \textit{reduced} HOMFLY homology, obtained by killing the action of a certain polynomial algebra naturally associated to points on a link. Our first major result gives an invariant construction of reduced colored HOMFLY homology for arbitrary links.

\begin{thm} \label{thm: reduced_homology_links}
Given any colored link $\L = \L_1 \sqcup \dots \sqcup \L_d$ with components labeled $k_1, \dots, k_d$ and any point $p \in \L_i$ on any component, there exists a \textit{reduced} colored HOMFLY homology $\overline{HHH}(\L(p))$ related to the ordinary (unreduced) colored HOMFLY homology $HHH(\L)$ by

\begin{equation} \label{eq: conj1}
HHH(\L) \cong HHH(\bigcirc_{k_i}) \otimes \overline{HHH}(\L(p)).
\end{equation}

\end{thm}

Theorem \ref{thm: reduced_homology_links} appears below as Theorems \ref{thm: link_invariance} and \ref{thm: main_reduction_thm}. Partial progress had been made by Wedrich in \cite{Wed19}, who constructed an invariant reduced version of colored HOMFLY homology for non-split links reduced at a component of minimal color (but did not prove \eqref{eq: conj1}). Our construction agrees with Wedrich's under these restrictions. Our method of proof mostly agrees with that of Wedrich, translated from the language of matrix factorizations to that of derived categories, together with a nonstandard set of generators and relations for singular Soergel bimodules and a study of the delicate interplay between derived and homotopy categories along different gradings.

To state our next result, we record the dimension of the triply-graded vector space $\overline{HHH}(\L(p))$ as a Laurent series in three variables $A, Q, T$. For experts, we note that these are the \textit{geometric} Hochschild, quantum, and homological gradings, related to the \textit{algebraic} gradings (in lowercase) by $A = aq^{-2}$, $Q = q^2$, $T = t^2q^{-2}$. We prove a certain \textit{refined exponential growth} property for the $k$-colored positive $(m, n)$ torus knots $T(m, n)_k$, originally conjectured in \cite{GGS18}.

\begin{thm} \label{thm: exp_grow}
The reduced HOMFLY homology of $T(m,n)_k$ satisfies

\begin{equation} \label{eq: conj2}
\dim \left( \overline{HHH} \left( T(m, n)_k(p) \right) \right) \Big|_{Q = 1} = \left( \dim \left( \overline{HHH} \left( T(m, n)_1(p) \right) \right) \Big|_{Q = 1} \right)^k
\end{equation}

\end{thm}

Theorem \ref{thm: exp_grow} appears below as Theorem \ref{thm: doubly_graded_growth}. This refines another result of Wedrich in \cite{Wed19}, who established a general bound on the ungraded total dimension of the form

\[
\dim \left( \overline{HHH} (\mathcal{K}_k(p)) \right) \Big|_{A = Q = T = 1} \geq \left( \dim \left( \overline{HHH} (\mathcal{K}_1(p)) \right) \Big|_{A = Q = T = 1} \right)^k
\]
for arbitrary knots $\mathcal{K}$. GGS also predict a triply-graded, and even a quadruply-graded, refinement of \eqref{eq: conj2}, relying on a modified $Q$-grading and a conjectural splitting of the $T$-grading into two components. In earlier work \cite{Wed16}, Wedrich also proved a triply-graded analog of \eqref{eq: conj2} for 2-component rational tangles at the level of Poincare polynomials with respect to a certain generating set of objects.  We do not address these extensions here, though our method of proof could conceivably be pushed to a triply-graded version. We discuss this method in Section \ref{subsec: computations} below.

The third and final conjecture which we address here concerns the behavior of $\overline{HHH}(\L(p))$ when $p$ lies on an unknotted component $\L_i$ of $\L$. GNSSS conjecture that the reduced homology in this setting splits into \textit{homological blocks}, independently of the color $k_i$ of $\L_i$, which satisfy certain \textit{color shifting} properties as $k_i$ changes.

\begin{conj}[\cite{GNSSS16}] \label{conj: color_shift}

Let $\L = \L_1 \sqcup \dots \sqcup \L_d$ be a link with $\L_1 \sim \bigcirc$ and suppose $p \in \L_1$. For any choice of colors $k_1, \dots, k_d$, denote the graded dimension of $\overline{HHH}$ reduced at $p$ of the corresponding colored link by $\mathcal{P}_{k_1, \dots, k_d}(A, Q, T)$. Then there exists some polynomial $f_{k_2, \dots, k_d}(A, Q, T, x)$, independent of $k_1$, which satisfies

\begin{equation} \label{eq: conj3}
\mathcal{P}_{k_1, \dots, k_d}(A, Q, T) = f_{k_2, \dots, k_d}(A, Q, T, x = Q^{k_1})
\end{equation}
for all $k_1 \geq 0$ (where here $\mathcal{P}_{0, k_2, \dots, k_d} := \dim(HHH(\L - \L_1))$).

\end{conj}

At the decategorified level of colored HOMFLY polynomials, Conjecture \ref{conj: color_shift} was proven by Wedrich in \cite{Wed16} by first computing the image of $\L - \L_1$ in the HOMFLY skein of the annulus, then pairing with the unknot component $\L_1$. While this technique is not available at the categorified level, we do explicitly verify \eqref{eq: conj3} for positive torus links $T(2, 2n)$ and $T(3, 3)$ over a broad range of colors, up to $k_1 = 90$ in many cases. This is the first confirmation of Conjecture \ref{conj: color_shift} for any nontrivial links; indeed, even evaluating the left-hand side of \eqref{eq: conj3} requires a version of Theorem \ref{thm: reduced_homology_links} for reduction at non-minimal color.

\subsection{Computations} \label{subsec: computations}

We establish Theorem \ref{thm: exp_grow} by explicitly computing the reduced homology of column-colored positive torus knots. Using Theorem \ref{thm: reduced_homology_links}, it suffices to compute the unreduced homology; this was accomplished in our previous work \cite{Con24}. There, we use a similar comparison theorem to reduce to the computation of an even larger invariant, the so-called ``$y$-ified column projector-colored HOMFLY homology". We compute this latter invariant directly in \cite{Con23}, generalizing previous work of Elias--Hogancamp--Mellit in \cite{EH19, HM19}.

The result is an explicit recursive formula for $\dim \left( \overline{HHH}(T(m, n)_k(p)) \right)$; see Theorem \ref{thm: red_torus_links} for the exact statement. In fact, our formula suffices to compute the colored HOMFLY homology of all positive torus \textit{links} with at most one nontrivially-colored component. This formula is well-adapted for evaluation with a computer, and we have implemented a program in SageMath which computes these invariants \cite{Con24b}. We can carry out huge computations using this program; for example, we can compute the homology of $3$-colored $T(4, 7)$ (which has total dimension $1,685,159$) and $1$-colored $T(13, 14)$ (dimension: $71,039,373$) in a few seconds. This program is freely available from the author's website, and we encourage further exploration!

All of our evidence for Conjecture \ref{conj: color_shift} comes directly from these computations. On the other hand, Theorem \ref{thm: exp_grow} comes from a careful analysis of the recursion of Theorem \ref{thm: red_torus_links}. After setting $Q = 1$ and taking some care to deal with denominators of the form $(1 - Q)^{-1}$, this recursion dramatically simplifies; see Proposition \ref{prop: simpler_recursion}. We prove directly that the resulting recursion exhibits exponential growth when $\gcd(m, n) = 1$. To do this, we once again use in an essential way that our recursion computes the much larger projector-colored homology. This latter invariant implements coloring in part by cabling, and we show that the recursion restricts in a precise sense to each link component of the corresponding cable.

\subsection{Future Directions}

Many of the GGS and GNSSS conjectures remain outstanding. We list some of these conjectures related to the work in this paper in no particular order.

\subsubsection{Color Shifting}

Conjecture \ref{conj: color_shift} remains wide open in full generality. In particular, it would be extremely interesting to investigate whether this behavior holds for non-torus links, as computing even the trivially-colored HOMFLY homology of these remains prohibitively difficult. GNSSS also conjecture that all of their differentials act within each homological block, dramatically simplifying the resulting spectral sequences.

Already for torus links, this conjecture could have dramatic consequences. Suppose one knew that the color shifting conjecture were true. By computing the colored homology of $T(n, dn)$ with a single nontrivially colored component for some small colors, one could then deduce a general formula for each homological block. The total homology is known to stabilize as $d \to \infty$, and one could imagine that each block would stabilize too. If it were known that the $\mathfrak{gl}(N)$ differentials $d_N$ respected each homological block, we would have additional restrictions on which terms could possibly cancel to produce the stable $\mathfrak{gl}(N)$ homology of torus knots. Computing this homology remains a major outstanding problem; see e.g. the recent work of Ballinger--Gorsky--Hogancamp--Wang \cite{BGHW25}.

\subsubsection{Reduced Homology for Arbitrary Colors}

This paper deals with \textit{column}-colored homology, in which link components carry positive integers $k$ which correspond to partitions $(1^k)$ in the stable representation theory of $\mathfrak{sl}(N)$. To our knowledge, there is still no construction of reduced colored HOMFLY homology for other partitions. Most of the constructions of unreduced colored HOMFLY for these colors pass through the categorified projectors of Cautis \cite{Cau17} and Elias--Hogancamp \cite{EH17b}; see the introduction of \cite{Con24} for a catalog of these theories. While we deal with ``intrinsically" column-colored homology here, there is also a ``projector-colored" version of this invariant, and Theorem \ref{thm: reduced_homology_links} admits a modification relating $\overline{HHH}$ to this version. We would expect a similar relationship in general. This would be especially interesting for non-rectangular colors, as GGS and GNSSS typically restrict to the rectangular case.

\subsubsection{Mirror Symmetry}

GGS conjecture that the reduced colored HOMFLY homologies for a knot colored by a partition $\lambda$ and its transpose $\lambda^t$ should be related by exchanging the $Q$ and $T$ gradings. This would categorify a well-known symmetry of the colored HOMFLY polynomial due to Liu and Peng \cite{LP09, LP11}. This behavior is surprising, as the $Q$ and $T$ gradings play dramatically different roles in the construction of HOMFLY homology. For uncolored HOMFLY homology (in which $\lambda = \square$ is stable under transposition), this was resolved for knots by Oblomkov--Rozansky and for links by Gorsky--Hogancamp--Mellit, as discussed in Section \ref{subsec: intro}. In the colored setting, mirror symmetry was shown to hold for row- and column-colored homology of positive torus knots in previous work of the author \cite{Con23}, but nothing is known outside of this setting. It would be especially interesting to see whether this result holds outside of rectangular colors.

\subsubsection{Supergroup Differentials}

The DGR conjectures included a $\mathbb{Z}$-indexed collection of differentials $d_N$ acting on uncolored HOMFLY homology, inducing spectral sequences as described in Section \ref{subsec: intro}. Among the GGS conjectures is a refinement of this collection to a $\mathbb{Z} \times \mathbb{Z}$-indexed collection of differentials $d_{m|n}$ acting on colored HOMFLY homology. These differentials are expected to again induce spectral sequences to $\mathfrak{gl}(m|n)$ homology, a family of conjectural invariants categorifying the $\mathfrak{gl}(m|n)$ Reshetikhin--Turaev link polynomials. Even understanding $\mathfrak{gl}(1|1)$ homology within this framework would be extremely interesting, as this is widely expected to agree with knot Floer homology. See \cite{BPRW22} for recent results in this direction.

\subsubsection{Quadruply-Graded Homology}

The titular conjecture of \cite{GGS18} is the existence of a $4^{th}$ grading on colored HOMFLY homology, refining the homological $t$-grading into distinct $t_c$ ($c$ for ``column") and $t_r$ ($r$ for ``row") gradings. Nearly all of the GGS and GNSSS conjectures admit cleaner or more general quadruply-graded versions, so such a construction would open up a world of possibilities. Nevertheless, this conjecture remains wide open.

\subsection{Organization}

In Section \ref{sec: background}, we recall the necessary homological algebra for constructing colored HOMFLY homology. Most of this material is well-known; the most interesting result is Corollary \ref{cor: cone_a_for_cone_t} regarding the interplay between the derived and homotopy categories in $K^b(D^b(\SSBim))$. We construct $\overline{HHH}$ and prove Theorem \ref{thm: reduced_homology_links} in Section \ref{sec: link_hom}. As mentioned above, we rely on a nonstandard presentation of singular Soergel bimodules via generators and relations which may be of independent interest; see Definition \ref{def: new_web_rels} for details. Finally, in Section \ref{sec: torus} we specialize to the case of positive torus links, proving Theorem \ref{thm: exp_grow} and exhibiting some evidence for Conjecture \ref{conj: color_shift}.

\subsection{Acknowledgements}

We would like to thank William Ballinger, Anna Beliakova, Elijah Bodish, Ben Elias, Eugene Gorsky, Matt Hogancamp, Lukas Lewark, David Rose, Paul Wedrich, Arno Wildi, and Haihan Wu for useful discussions. Special thanks to Josh Wang, who originally brought Theorem \ref{thm: exp_grow} to our attention and suggested that we might implement Theorem \ref{thm: red_torus_links} on a computer.

During the preparation of this work, the author was partially supported by NSF CAREER Grant DMS-2144463 and the Simons Collaboration ``New Structures in Low-Dimensional Topology".

\subsection{Conventions}

There are various overall normalization conventions in the literature which shift the value of $HHH(\L)$ by some monomial $a^it^jq^k$. Because we rely heavily on computer algebra involving rational functions, it quickly becomes very difficult to track these absolute shifts when comparing different invariants. With this in mind, all statements involving link invariants in this paper should be presumed to hold up to some relative shift of this form.

In Section \ref{sec: torus}, we make frequent use of the \textit{geometric} variables $A := aq^{-2}$, $Q := q^{2}$, $T := t^2q^{-2}$. In Theorem \ref{thm: red_torus_links}, we also make use of the \textit{quantum factorial} $[k]! = [k] [k-1] \dots [2] [1]$, where $[j] := q^{j - 1} + q^{j - 3} + \dots + q^{-j + 3} + q^{-j + 1}$.

\section{Background} \label{sec: background}

\subsection{Koszul Complexes}

In this section we recall some results in homological algebra. Most of this material is quite standard and can be skipped by experts. The primary novelty is Corollary \ref{cor: cone_a_for_cone_t}, in which we prove a relationship between Koszul complexes in two distinct homological directions.

\begin{notn}
Given a graded ring $R$, denote by $R-\Mod{}$ the category of graded $R$-modules and (degree $0$) graded $R$-module homomorphisms. We denote the grading group on $R$ by $\mathbb{Z}_q$ and write the degree of elements multiplicatively using powers of $q$. We also use powers of $q$ to represent grading shift functors on $R-\Mod{}$. For example, given a module $M \in R-\Mod{}$, the module $q^3M$ has graded components $(q^3M)_i = M_{i - 3}$. We will often call this grading the \textit{quantum grading} or the \textit{internal grading}.

We denote by $D^b(R-\Mod{})$ the bounded derived category of graded $R$-modules. We denote the derived homological grading group by $\mathbb{Z}_a$, and we again write degrees and grading shift functors multiplicatively using powers of $a$. We will often call this grading the \textit{Hochschild grading} or the \textit{derived grading}. We will use the notation $C \cong_{q.i.} C'$ for isomorphism in $D^b(R-\Mod{})$ (i.e. quasi-isomorphism).

Finally, we denote by $K^b(D^b(R-\Mod{}))$ the bounded homotopy category of $D^b(R-\Mod{})$. This category has a new homological grading group which we denote by $\mathbb{Z}_t$. We will similarly write degrees of elements and grading shift functors using powers of $t$. We call the $\mathbb{Z}_t$ grading the \textit{homological grading}. We will use the notation $C \simeq C'$ for isomorphism in $K^b(D^b(R-\Mod{}))$ (i.e. homotopy equivalence).
\end{notn}

\begin{defn} \label{def: Koszul_complex}
Fix a graded ring $R$, a complex $C \in D^b(R-\Mod{})$ and a sequence of homogeneous elements $(r_1, \dots, r_n) \in R$ of degrees $\deg(x_i) = q^{d_i}$. For each $1 \leq i \leq n$, let $\theta_i$ be a formal variable of degree $a^{-1}q^{d_i}$. Then the \textit{Hochschild Koszul complex} on $C$ for the elements $r_1, \dots, r_n$ is the complex

\[
K_a(C; r_1, \dots, r_n) := \big(C \otimes \bigwedge[\theta_1, \dots, \theta_n], d_C \otimes 1 + \sum_{i = 1}^n r_i \otimes \theta_i^{\vee} \big) \in D^b(R-\Mod{}).
\]

Now fix a complex $C \in K^b(D^b(R-\Mod{}))$. For each $1 \leq i \leq n$, let $\xi_i$ be a formal variable of degree $t^{-1}q^{d_i}$. Then the \textit{homological Koszul complex} on $C$ for the elements $r_1, \dots, r_n$ is the complex

\[
K_t(C; r_1, \dots, r_n) := \big(C \otimes \bigwedge[\xi_1, \dots, \xi_n], d_C \otimes 1 + \sum_{i = 1}^n r_i \otimes \theta_i^{\vee} \big) \in K^b(D^b(R-\Mod{})).
\]
\end{defn}

\begin{rem}
It's worth pointing out that Koszul complexes respect quasi-isomorphism and homotopy equivalences, so Definition \ref{def: Koszul_complex} extends to produce functors on the categories $D^b(R-\Mod{})$ and $K^b(D^b(R-\Mod{}))$. For example, suppose we have  a morphism from $C$ to $C'$ in $D^b(R-\Mod{})$. Then there is a roof of the form

\begin{center}
\begin{tikzcd}
& D \arrow[dl, "f"'] \arrow[dr, "g"] & \\
C & & C'
\end{tikzcd}
\end{center}
where $f$ is a quasi-isomorphism. Applying $K_a(-; r)$ for an element $r \in R$ of degree $d$, we obtain a roof of the form

\begin{center}
\begin{tikzcd}
& \theta D \arrow[r, "r"] \arrow[dl, "-f"'] \arrow[dr, "-g",pos=.9] & D \arrow[dl, "f"',pos=.9] \arrow[dr, "g"] & \\
\theta C \arrow[r, "r"] & C & \theta C' \arrow[r, "r"] & C'
\end{tikzcd}
\end{center}

For this roof to represent a morphism from $K_a(C; r)$ to $K_a(C'; r)$ in $D^b(R-\Mod{})$, we need to check that the diagonal leftwards arrows constitute a quasi-isomorphism from $K_a(D; r)$ to $K_a(C; r)$. These arrows certainly form a chain map $K_a(f)$ (since $f$ is $R$-linear), and $K_a(f)$ is a quasi-isomorphism if and only if $\mathrm{Cone}(K_a(f))$ is acyclic. We can write $\mathrm{Cone}(K_a(f))$ as the complex

\begin{center}
\begin{tikzcd}
a^{-1} \theta D \arrow[r, "r"] \arrow[d, "-f"'] & a^{-1} D \arrow[d, "f"] \\
\theta C \arrow[r, "r"] & C 
\end{tikzcd}
\end{center}

This is exactly the mapping cone for the action of $r$ on the mapping cone $\mathrm{Cone}(f)$, so it suffices to show that this action is a quasi-isomorphism. But since $f$ is a quasi-isomorphism, $\mathrm{Cone}(f)$ is itself acyclic, so any endomorphism on $\mathrm{Cone}(f)$ is a quasi-isomorphism.
\end{rem}

Koszul complexes are particularly useful for constructing projective resolutions when the sequence $(r_1, \dots, r_n)$ is regular, as we recall next.

\begin{defn} \label{def: reg_seq}
A sequence $(r_1, \dots, r_n) \subset R$ is called \textit{regular} if $r_i$ is not a zero divisor in $R/\langle r_1, \dots, r_{i - 1} \rangle$ for each $1 \leq i \leq n$.
\end{defn}

\begin{prop} \label{prop: simplifying_koszul}
Fix a commutative ring $R$ and a regular sequence $\{r_1, \dots, r_n\} \subset R$, and set $S := R/\langle r_1, \dots, r_n \rangle$. Then for every complex $C \in D^b(R-\Mod{})$, there is a quasi-isomorphism

\[
K_a(C; r_1, \dots, r_n) \cong_{qi} C \tensor[R, der] S \in D^b(R-\Mod{}).
\]
\end{prop}

\begin{proof}
This is more or less the definition of $\tensor[R, der]$. To see this, recall that $C \tensor[R, der] S$ can be computed by first passing to a projective resolution of the $R$-module $S$. By standard homological algebra, we can take the Koszul resolution $K_a(R; r_1, \dots, r_n) \cong_{qi} S \in D^b(R-\Mod{})$. Taking the tensor product with $C$ gives exactly $K_a(C; r_1, \dots, r_n)$.
\end{proof}

In general, whether a sequence $(r_1, \dots, r_n)$ is regular might depend on the order of the elements. Fortunately, we will never have to deal with this ambiguity due to the following standard result on graded rings; see e.g. Theorem 27 of \cite{Mat80} and the subsequent remarks.

\begin{lem} \label{lem: reg_seq_perm}
Let $R = \bigoplus_{i \geq 0} R_i$ be a non-negatively graded ring, and suppose $(r_1, \dots, r_n)$ is a regular sequence in $R$ consisting of homogeneous elements of positive degree. Then any permutation of $(r_1, \dots, r_n)$ is a regular sequence in $R$.
\end{lem}

Next, we introduce the main technical tools we use to manipulate Koszul complexes. Many of the statements below hold in both the Hochschild and the homological setting. In this case, we do not distinguish between the two, using for example $K(C; r)$ to denote either $K_a(C; r)$ for a complex $C \in D^b(R-\Mod{})$ or $K_t(C; r)$ for a complex $C \in K^b(D^b(R-\Mod{}))$.

\begin{lem} \label{lem: kosz_homotopic}
Fix a complex $C$, and suppose the actions of $r$ and $r'$ on $C$ are homotopic via some homotopy $h$ satisfying $[d_C, h] = r - r'$. Then the map

\[
\Psi := 1 \otimes 1 + h \otimes \theta^{\vee} \colon K(C; r) \to K(C; r')
\]
is a chain isomorphism.
\end{lem}

\begin{proof}
Direct computation.
\end{proof}

\begin{cor} \label{cor: row_ops}
Suppose $s_1, s_2, r \in R$ are homogeneous elements satsifying $\mathrm{deg}(rs_1) = \mathrm{deg}(s_2)$. Then for any complex $C$, we have a chain isomorphism $K(C; s_1, s_2) \cong K(C; s_1, rs_1 + s_2)$.
\end{cor}

\begin{proof}
By definition, we have a chain isomorphism $K(C; s_1, s_2) \cong K(K(C; s_1); s_2)$. The actions of $s_2$ and $rs_1 + s_2$ on $K(C; s_1)$ are nullhomotopic via the homotopy $r \theta_1$, which satisfies $[d, r \theta_1] = rs_1$. The desired isomorphism then follows immediately from Lemma \ref{lem: kosz_homotopic}.
\end{proof}

We close with some useful tools for relating Koszul complexes in the homological and Hochschild directions.

\begin{prop} \label{prop: quot_for_cone}
Fix a commutative ring $S$, and set $R = S[x]$. Then for each complex $C \in K^b(D^b(R-\Mod{}))$, we have

\[
C \tensor[R,der] S \simeq \Cone_t(x) \in K^b(D^b(S-\Mod{})).
\]
Here $S$ is considered as an $R$-module via extending scalars along the quotient map $\pi \colon R \to R/\langle x \rangle \cong S$, and we implicitly include $K^b(D^b(R-\Mod{})) \hookrightarrow K^b(D^b(S-\Mod{}))$ via restricting scalars along the inclusion map $\iota \colon S \hookrightarrow R$.
\end{prop}

\begin{proof}
We adapt the proofs of analogous propositions provided in \cite{Wed19}, Lemma 3.39 and \cite{Ras15}, Lemma 5.15 from the setting of matrix factorizations to the derived category.

Pick a resolution (at the derived level $D^b(R-\Mod{})$) $C^f$ of $C$ by free $R$-modules; then we can compute $C \tensor[R,der] S$ via the usual underived tensor product $C^f \otimes_R S$. We also model $\Cone_t(x)$ using this resolution, i.e.

\[
\Cone_t(x) := t^{-1} C^f \xrightarrow{\quad x \quad} C^f.
\]

Now since every chain module of $C^f$ is free over $R$, $C^f$ admits a decreasing filtration $C^f = C^f_0 \supset C^f_1 \supset \dots$ by $x$-degree, and $C^f \otimes_R S$ is exactly the quotient $C^f_0/C^f_1$ with respect to this filtration. Multiplication by $x$ maps $C^f_0$ onto $C^f_1$, and this map is an $S$-linear isomorphism. In particular, we can Gaussian eliminate along this component of the differential of $\Cone_t(x)$ in $K^b(D^b(S-\Mod{}))$, and what remains is exactly the (unshifted) quotient $C^f_0/C^f_1$.
\end{proof}

Combining Propositions \ref{prop: simplifying_koszul} and \ref{prop: quot_for_cone} gives a useful consequence.

\begin{cor} \label{cor: cone_a_for_cone_t}
Fix a commutative ring $S$, and set $R = S[x_1, \dots, x_n]$. Then for each complex $C \in K^b(D^b(R-\Mod{}))$, we have

\[
K_a(C; x_1, \dots, x_n) \simeq K_t(C; x_1, \dots, x_n) \in K^b(D^b(S-\Mod{})).
\]
\end{cor}

\begin{proof}
Since $(x_1, \dots, x_n)$ is a regular sequence in $R$, we can apply Proposition \ref{prop: simplifying_koszul} to obtain a quasi-isomorphism $K_a(C^i; x_1, \dots, x_n) \cong_{qi} C^i \tensor[R, der] S$ in $D^b(S-\Mod{})$ for each complex $C^i \in D^b(R-\Mod{})$ in homological degree $i$. These quasi-isomorphisms assemble to produce a homotopy equivalence $K_a(C; x_1, \dots, x_n) \simeq C \tensor[R, der] S$ in $K^b(D^b(S-\Mod{}))$. The result then follows directly from Proposition \ref{prop: quot_for_cone}.
\end{proof}

\subsection{Symmetric Polynomials}

\begin{defn} \label{def: composition}
Let $N$ be a nonnegative integer. A \textit{composition} of $N$ is a (possibly empty) ordered sequence $\aa = (a_1, \dots, a_n)$ of positive integers satisfying $\sum_{i = 1}^n a_i = N$. We write $\aa \vdash N$ to indicate that $\aa$ is a composition of $N$.
\end{defn}

\begin{defn} \label{def: symm_polys}
Let $\X$ be an alphabet of formal variables $\X = \{x_1, \dots, x_N\}$, and consider the polynomial algebra $R[\X]$ as a $\mathbb{Z}_q$-graded ring by setting $\deg(x_i) = q^2$. The symmetric group $\SG_N$ naturally acts on $R[\X]$ by permuting the variables. Given a composition $\aa = (a_1, \dots, a_n) \vdash N$, we denote by $\Sym^{\aa}(\X)$ the invariant subalgebra of $R[\X]$ under the action of the parabolic subgroup $\SG_{\aa} := \SG_{a_1} \times \dots \times \SG_{a_n}$ of $\SG_N$.
\end{defn}

\begin{ex}
When $\aa = (1, 1, \dots, 1) \vdash N$, the parabolic subgroup $\SG_{\aa}$ is trivial, and $\Sym^{\aa}(\X) = R[\X]$.
\end{ex}

\begin{ex}
When $\aa = (N) \vdash N$, $\SG_{\aa} = \SG_N$, and $\Sym^{\aa}(\X)$ is the ring of symmetric polynomials in $\X$, which we denote by $\Sym(\X)$.
\end{ex}

We'll make use of a few distinguished generating sets of $\Sym(\X)$.

\begin{defn} \label{def: elem_symms}
For each $i \geq 0$, the \textit{$i^{th}$ elementary symmetric polynomial} $e_i(\X) \in \Sym(\X)$ is the sum of all square-free monomials in the alphabet $\X$ of polynomial degree $i$. The \textit{$i^{th}$ complete symmetric polynomial} $h_i(\X) \in \Sym(\X)$ is the sum of all monomials in the alphabet $\X$ of polynomial degree $i$.
\end{defn}

\begin{prop}[Fundamental Theorem of Symmetric Polynomials] \label{prop: fund_thm_symm}
Let $\X$ be an alphabet of size $N$. Then $\Sym(\A) \cong R[e_1(\X), \dots, e_N(X)] \cong R[h_1(\X), \dots, h_N(\X)]$.
\end{prop}

For an alphabet $\X$ of size $N$, $e_i(\X) = 0$ for all $i > N$. We will nevertheless find it useful to artificially extend $\Sym(\X)$ to include these elements using the following result.

\begin{cor} \label{cor: artificial_size}
Let $\A$, $\B$ be two alphabets of size $|\A| = m \leq n = |\B|$. Set 
\[
\IS := \langle e_i(\B) \ | \ m < i \leq n \rangle \subset \Sym(\B).
\]
Then there is an algebra isomorphism $\Sym(\A) \cong \Sym(\B)/\IS$ taking $e_i(\A)$ to (the equivalence class of) $e_i(\B)$.
\end{cor}

Given an arbitrary composition $\aa = (a_1, \dots, a_n) \vdash N$, we have an isomorphism of graded algebras

\[
\Sym(\X_1) \otimes \dots \otimes \Sym(\X_n)\cong \Sym^\aa(\X) ,
\]
where $\X_i$ is an alphabet of size $a_i$ for each $i$. Considering $\X$ as a union of the alphabets $\X_i$, this isomorphism sends a pure tensor of the form $e_{k_1}(\X_1) \otimes \dots \otimes e_{k_n}(\X_n)$ to the product $\prod_{i = 1}^n e_{k_i}(\X_i)$. To emphasize this decomposition, we will often use the notation

\[
\Sym(\X_1 | \X_2 | \dots | \X_n) := \Sym^{\aa}(\X).
\]

More generally, given two (not necessarily disjoint) alphabets $\A$ and $\B$ of sizes $m$ and $n$, respectively, we will often write $\Sym(\A | \B)$ to denote the partially symmetric polynomial ring $\Sym^{(m, n)}(\A \cup \B)$. In this ring, beyond the usual elementary and complete symmetric polynomials in either alphabet, there are also very useful linear combinations of such, typically referred to as elementary symmetric polynomials in linear combinations of alphabets. We recall here only the small portion of this extremely useful formalism which we will need; see e.g. \S 2 of \cite{HRW24} for a more thorough exposition.

\begin{defn}
Given two (not necessarily disjoint) alphabets $\A, \B$, for each $i \geq 0$, we set

\begin{align*}
e_i(\A + \B) & := \sum_{j + k = i} e_j(\A) e_k(\B) \in \Sym(\A | \B); \\
e_i(\A - \B) & := \sum_{j + k = i} (-1)^k e_j(\A) h_k(\B) \in \Sym(\A | \B).
\end{align*}
\end{defn}

We will make frequent use of the following cancellation property:

\begin{prop} \label{prop: solving_for_alphabets}
Let $\A$, $\B$, $\C$ be three (not necessarily disjoint) alphabets. Then for each $k \geq 0$, we have

\[
\langle e_i(\A - \C) - e_i(\B - \C) | 1 \leq i \leq k \rangle = \langle e_i(\A) - e_i(\B) | 1 \leq i \leq k \rangle \in \Sym(\A | \B | \C).
\]
\end{prop}

\begin{proof}
By induction, it suffices to show

\[
\left((e_i(\A - \C) - e_i(\B - \C)\right) - \left(e_i(\A) - e_i(\B) \right) \in \langle e_j(\A) - e_j(\B) | 1 \leq j < i \rangle
\]
for each $i \geq 1$. We proceed by direct computation:

\begin{align*}
(e_i(\A - \C) - e_i(\B - \C)) & = \left( \sum_{j = 0}^i (-1)^{i - j} e_j(\A) h_{i-j}(\C) \right) - \left( \sum_{j = 0}^i (-1)^{i - j} e_j(\B) h_{i-j}(\C) \right) \\
& = (e_i(\A) - e_i(\C)) + \sum_{j = 0}^{i - 1} (-1)^{i - j} h_{i-j}(\C) \left(e_j(\A) - e_j(\B) \right).
\end{align*}
\end{proof}

\begin{cor} \label{cor: solving_for_alphabets}
Let $\A$, $\B$, $\C$ be three (not necessarily disjoint) alphabets. Then for each $k \geq 0$, we have

\[
\langle e_i(\A) - e_i(\B + \C) | 1 \leq i \leq k \rangle = \langle e_i(\A - \B) - e_i(\C) | 1 \leq i \leq k \rangle \in \Sym(\A | \B | \C).
\]
\end{cor}

\begin{proof}
By Proposition \ref{prop: solving_for_alphabets}, we immediately obtain

\[
\langle e_i(\A - \B) - e_i((\B + \C) - \B) | 1 \leq i \leq k \rangle = \langle e_i(\A) - e_i(\B + \C) | 1 \leq i \leq k \rangle.
\]

It therefore suffices to show that $e_i((\B + \C) - \B) = e_i(\C)$ for each $i$. Expanding the left-hand expression, we obtain

\[
e_i((\B + \C) - \B) = \sum_{j + k + l = i} e_j(\C) \left( (-1)^l e_k(\B) h_l(\B) \right).
\]

The desired result then follows from the well-known identity $\sum_{l + k = n} (-1)^l e_k(\B) h_l(\B) = 0$ for all $n \geq 1$.
\end{proof}

\section{Colored HOMFLY Homology} \label{sec: link_hom}

\subsection{Singular Soergel Bimodules} \label{subsec: SSBim}

\begin{defn} \label{def: Bim}
For each $N \geq 0$, we denote by $\Bim_N$ the $2$-category whose objects are compositions $\aa \vdash N$, $1$-morphisms from $\aa$ to $\bb$ are graded $(\Sym^\bb(\X), \Sym^\aa(\X))$-bimodules, and $2$-morphisms are graded bimodule homomorphisms. Horizontal composition of $1$- and $2$-morphisms (denoted $\star$) along objects $\aa \to \bb \to \cc$ is given by tensor product over $\Sym^\bb(\X)$, and vertical composition of $2$-morphisms is composition of homomorphisms. The $2$-categories $\Bim_N$ for various $N$ assemble into a monoidal $2$-category $\Bim := \bigsqcup_{N \geq 0} \Bim_N$ with tensor product 

\[
\boxtimes \colon \Bim_{N_1} \times \Bim_{N_2} \to \Bim_{N_1 + N_2}
\]
given on objects by concatenation (i.e. $\aa \boxtimes \bb = (\aa, \bb)$) and on $1$- and $2$-morphisms by tensor product over $R$.

We denote by $\Bim_\aa^\bb$ the ordinary category $\Hom_\Bim(\aa, \bb)$.

\end{defn}

When $\SG_{\aa} \subset \SG_\bb$, we can view $\Sym^\aa(\X)$ as an object of $\Bim_\aa^\bb$ (resp. $\Bim_\bb^\aa$) with left action (resp. right action) given by restricting scalars along the inclusion $\Sym^\bb(\X) \subset \Sym^\aa(\X)$. In this case, we have distinguished \textit{merge} and \textit{split} bimodules

\[
_{\bb}M_{\aa} := q^{\ell(\bb) - \ell(\aa)} \mathrm{Sym}^{\bb}(\mathbb{X}) \otimes_{\mathrm{Sym}^{\bb}(\mathbb{X})} \mathrm{Sym}^{\aa}(\mathbb{X}); \quad _{\aa}S_{\bb} := \mathrm{Sym}^{\aa}(\mathbb{X}) \otimes_{\mathrm{Sym}^{\bb}(\mathbb{X})} \mathrm{Sym}^{\bb}(\mathbb{X}).
\]
Here $\ell(\aa)$ is the length of the longest word in the parabolic subgroup $\SG_\aa \subset \SG_N$ for each composition $\aa \vdash N$.

\begin{defn} \label{def: SSBim}
A \textit{singular Bott-Samelson bimodule} is a bimodule which is generated by those of the form $_\bb M_\aa$ and $_\aa S_\bb$ under grading shifts, horizontal and vertical composition, and direct sums. The $2$-category $\SSBim$ of \textit{singular Soergel bimodules} is the smallest full monoidal $2$-subcategory of $\Bim$ containing all Bott-Samelson bimodules which is closed under taking direct summands. As with $\Bim$, we denote by $\SSBim_\aa^\bb$ the ordinary category $\Hom_\SSBim(\aa, \bb)$.
\end{defn}

Singular Bott-Samelson bimodules are most conveniently described using the diagrammatic language of Type A webs, as developed for instance by Cautis--Kamnitzer--Morrison \cite{CKM14}. In this description, merge and split bimodules for a pair of compositions $(a, b), (a + b) \vdash a + b$ are depicted using trivalent vertices as follows:

\begin{gather*}
_{(a + b)}M_{(a, b)} := 
\begin{tikzpicture}[anchorbase]
\draw[webs] (0,0) node[below]{$a$} to[out=90,in=180] (.5,.5);
\draw[webs] (1,0) node[below]{$b$} to[out=90,in=0] (.5,.5);
\draw[webs] (.5,.5) to (.5,1) node[above]{$a+b$};
\end{tikzpicture}
; \quad \quad _{(a, b)}S_{(a + b)} :=
\begin{tikzpicture}[anchorbase]
\draw[webs] (.5,0) node[below]{$a+b$} to (.5,.5);
\draw[webs] (.5,.5) to[out=180,in=270] (0,1) node[above]{$a$};
\draw[webs] (.5,.5) to[out=0,in=270] (1,1) node[above]{$b$};
\end{tikzpicture}
\end{gather*}

This description is extended to all singular Bott-Samelson bimodules by depicting horizontal composition of bimodules (i.e. $\star$) as vertical concatenation and the monoidal product $\boxtimes$ as horizontal concatenation. Bimodule homomorphisms between singular Bott-Samelson bimodules have a compatible description in this calculus as foams between webs, giving an equivalence between the diagrammatic and algebraic $2$-categories of such; see e.g. \cite{HRW21b}, Appendix 1 for a complete dictionary. 

\begin{ex}
Given a diagrammatic description of a singular Bott-Samelson bimodule as a (direct sum of) web(s), one can easily read off the tensor factors by considering the labels on horizontal cross-sections. These indicate the size of each alphabet in which the relevant polynomial algebras are partially symmetric when read from left to right. For example, the web below depicts the tensor product indicated on the right (up to a grading shift):

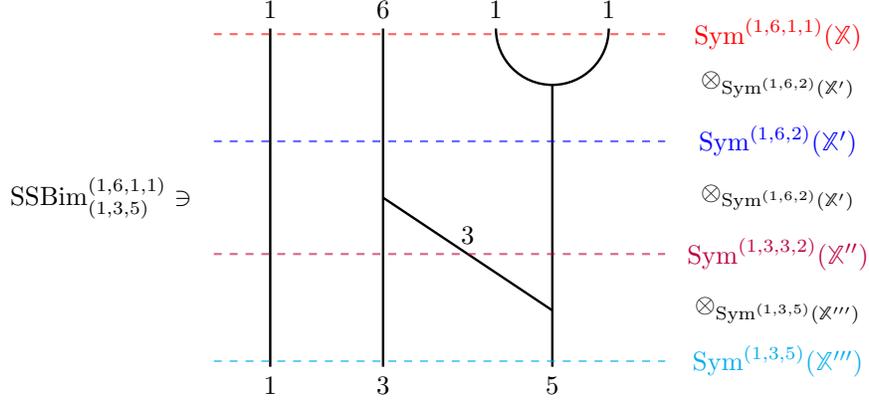
\begin{center}
\begin{tikzpicture}[anchorbase,scale=1.5]
	\node at (-1.5,1.5) {$\SSBim_{(1,3,5)}^{(1,6,1,1)} \ni$};
	\draw[webs] (0,0) node[below]{$1$} to (0,3) node[above]{$1$};
	\draw[webs] (1,0) node[below]{$3$} to (1,3) node[above]{$6$};
	\draw[webs] (2.5,0) node[below]{$5$} to (2.5,2.5);
	\draw[webs] (2.5,2.5) to[out=180,in=270] (2,3) node[above]{$1$};
	\draw[webs] (2.5,2.5) to[out=0,in=270] (3,3) node[above]{$1$};
	\draw[webs] (1,1.5) to node[midway,above]{$3$} (2.5,.5);
	
	\draw[thick, dashed, red, opacity=.6] (-.5,2.95) to (3.5,2.95);
	\draw[thick, dashed, blue, opacity=.6] (-.5,2) to (3.5,2);
	\draw[thick, dashed, purple, opacity=.6] (-.5,1) to (3.5,1);
	\draw[thick, dashed, cyan, opacity=.6] (-.5,.05) to (3.5,.05);
	
	\node at (4.5,2.95) [red]{$\mathrm{Sym}^{(1,6,1,1)}(\X)$};
	\node at (4.5,2.5) {$\otimes_{\mathrm{Sym}^{(1,6,2)}(\X')}$};
	\node at (4.5,2) [blue]{$\mathrm{Sym}^{(1,6,2)}(\X')$};
	\node at (4.5,1.5) {$\otimes_{\mathrm{Sym}^{(1,6,2)}(\X')}$};
	\node at (4.5,1) [purple]{$\mathrm{Sym}^{(1,3,3,2)}(\X'')$};
	\node at (4.5,.5) {$\otimes_{\mathrm{Sym}^{(1,3,5)}(\X''')}$};
	\node at (4.5,.05) [cyan]{$\mathrm{Sym}^{(1,3,5)}(\X''')$};
\end{tikzpicture}
\end{center}

Note that while we do not explicitly label each edge, the missing weights can be deduced easily using the ``Kirchhoff flow" condition that at each vertex, the sum of the incoming weights must equal the sum of the outgoing weights.
\end{ex}

Each trivalent vertex in a web indicates a tensor product over a partially symmetric polynomial algebra. We can also describe these tensor products by presenting a singular Bott-Samelson bimodule as an explicit quotient of a polynomial algebra in which we've identified those symmetric polynomials which slide through such vertices. The following presentation is well-known.

\begin{prop} \label{prop: web_quot_regular}
For each pair of positive integers $a, b$, let $\X$ be an alphabet of size $a$, $\Y$ an alphabet of size $b$, and $\M$ an alphabet of size $a + b$. Then we have algebra\footnote{These are not bimodule isomorphisms, as the actions on the left and right differ for merge and split bimodules (and are unspecified in the quotient).} isomorphisms

\[
q^{-ab} \ _{(a, b)} M_{(a + b)} \cong _{(a + b)}S_{(a, b)} \cong \cfrac{\Sym(\X | \Y | \M)}{\langle e_i(\X + \Y) = e_i (\M) | 1 \leq i \leq a + b \rangle}.
\]

Taking tensor products produces explicit descriptions of all singular Bott-Samelson bimodules as quotients of polynomial algebras. Moreover, the generators of the denominator in any such quotient form a regular sequence.
\end{prop}

\begin{proof}
The only thing that needs to be checked is that the generators in the denominator of the quotient form a regular sequence. By Lemma \ref{lem: reg_seq_perm}, this property is independent of their order. For a merge or split bimodule, each relation contains a unique summand of the form $e_i(\M)$. These terms are algebraically independent in $\Sym(\X | \Y | \M)$ by Proposition \ref{prop: fund_thm_symm}, so they certainly form a regular sequence. The result for a general singular Bott-Samelson bimodule follows by an easy induction on the number of vertices.
\end{proof}

\begin{defn} \label{def: ladder_web}
Fix $a, b \geq 0$. Then we set

\begin{gather*}
W_k := 
\begin{tikzpicture}[anchorbase]
\draw[webs] (0,0) node[below]{$a$}
   -- (0,2) node[pos=0.5,left] {}
   node[above] {$b$};
\draw[webs] (1,0) node[below]{$b$} -- (1,2) node[pos=.5, right]{$k$} node[above]{$a$};
\draw[webs] (1,.25) -- (0,.75) node[pos=.5,above]{};
\draw[webs] (0,1.25) -- (1,1.75) node[pos=.5,above]{};
\end{tikzpicture}
\in \SSBim_{(a,b)}^{(b,a)}.
\end{gather*}
We call singular Soergel bimodules of the form $W_k$ \textit{ladder webs}.
\end{defn}

We will find it convenient use a different presentation of ladder webs as quotients than the one given in Proposition \ref{prop: web_quot_regular}. This presentation is nonstandard, but it has the distinct advantage that all ladder webs $W_k$ for fixed $a$ and $b$ are quotients of a single polynomial algebra $\RS_{a,b}$. We will slightly abuse the diagrammatic notation here, presenting a ``web" with nonsensical labels which nevertheless captures the relations which we would like to impose.

\begin{notn}
For each pair of positive integers $a, b$, we assign the alphabets depicted on the left with sizes indicated on the right:

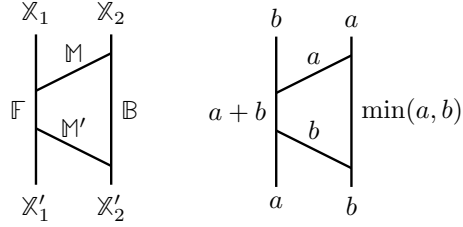
\begin{gather*}
\begin{tikzpicture}[anchorbase]
\draw[webs] (0,0) node[below]{$\X'_1$}
   -- (0,2) node[pos=0.5,left] {$\F$}
   node[above] {$\X_1$};
\draw[webs] (1,0) node[below]{$\X'_2$} -- (1,2) node[pos=.5, right]{$\B$} node[above]{$\X_2$};
\draw[webs] (1,.25) -- (0,.75) node[pos=.5,above]{$\M'$};
\draw[webs] (0,1.25) -- (1,1.75) node[pos=.5,above]{$\M$};
\end{tikzpicture}
\quad \quad
\begin{tikzpicture}[anchorbase]
\draw[webs] (0,0) node[below]{$a$}
   -- (0,2) node[pos=0.5,left] {$a + b$}
   node[above] {$b$};
\draw[webs] (1,0) node[below]{$b$} -- (1,2) node[pos=.5, right]{$\min(a, b)$} node[above]{$a$};
\draw[webs] (1,.25) -- (0,.75) node[pos=.5,above]{$b$};
\draw[webs] (0,1.25) -- (1,1.75) node[pos=.5,above]{$a$};
\end{tikzpicture}
\end{gather*}
\end{notn}

\begin{defn} \label{def: new_web_rels}
For each pair of positive integers $a, b$, set
 
\[
\RS_{a, b} := \Sym(\X_1 | \X_2 | \M | \M' | \F | \B | \X'_1 | \X'_2).
\]
We denote by $\IS_{a, b} \subset \RS_{a, b}$ the ideal of $\RS_{a, b}$ generated by the following elements:

\[
\IS_{a,b} \;=\;
\stretchleftright[1000]{<}
{\begin{aligned}
 & e_i(\X_2) - e_i(\M + \B) && 1 \le i \le a, \\
 & e_i(\X'_2) - e_i(\M' + \B) && 1 \le i \le b, \\
 & e_i(\F) - e_i(\X_1 + \M) && 1 \le i \le a+b, \\
 & e_i(\X_1 + \X_2) - e_i(\X_1' + \X_2') && 1 \le i \le \max(a,b)
\end{aligned}}
{>}
\]
Finally, for each $1 \leq k \leq \min(a, b)$, we let $\JS_k \subset \RS_{a, b}$ denote the ideal generated by the following elements:

\[
\JS_k \;=\;
\stretchleftright[1000]{<}
{\begin{aligned}
 & e_i(\M) && a-k < i \leq a, \\
 & e_i(\M') && b-k < i \leq b, \\
 & e_i(\B) && k < i \leq \min(a, b), \\
 & e_i(\X_1 + \X_2) - e_i(\X_1' + \X_2') && \max(a,b)+1 < i \leq a+b-k
\end{aligned}}
{>}
\]

\end{defn}

\begin{prop}
For any ladder web $W_k \in \SSBim_{(a, b)}^{(b, a)}$, there is an isomorphism\footnote{Up to an overall grading shift.} $W_k \cong \RS_{a, b}/(\IS_{a, b} \sqcup \JS_k)$ which preserves elementary symmetric polynomials in each alphabet.
\end{prop}

\begin{proof}
Consider the presentation of $W_k$ given in Proposition \ref{prop: web_quot_regular}. By Corollary \ref{cor: artificial_size}, we can artificially extend the size of the alphabets associated to the edges $\M, \M', \F$, and $\B$ to obtain the ring $\RS_{a,b}$ at the cost of imposing additional relations. This results in a presentation of $W_k$ as the quotient of $\RS_{a,b}$ by the following relations:

\begin{enumerate}
\item $e_i(\X_2) - e_i(\M + \B)$ \quad $1 \leq i \leq a$;

\item $e_i(\F) - e_i(\M + \X_1)$ \quad $1 \leq i \leq a+b-k$;

\item $e_i(\F) - e_i(\M' + \X'_1)$ \quad $1 \leq i \leq a+b-k$;

\item $e_i(\X'_2) - e_i(\M' + \B)$ \quad $1 \leq i \leq b$;

\item $e_i(\M)$ \quad $a-k < i \leq a$;

\item $e_i(\M')$ \quad $b-k < i \leq b$;

\item $e_i(\F)$ \quad $a+b-k < i \leq a+b$;

\item $e_i(\B)$ \quad $k < i \leq \min(a,b)$.
\end{enumerate}

Relations (1), (2), and (4) are present in $\IS_{a, b}$, while relations (5), (6), and (8) are present in $\JS_k$. For $i > a+b-k$, $e_i(\M + \X_1)$ vanishes by relation (5), so the excess relations $e_i(\F) - e_i(\X_1 + \M)$ with $i > a+b-k$ are equivalent to the relations (7).

This leaves the relations (3) from $W_k$ and the relations $\{e_i(\X_1 + \X_2) - e_i(\X_1' + \X_2') \ | \ 1 \leq i \leq a+b-k\}$ from $\IS_{a,b}$ and $\JS_k$. Using (2), we can replace (3) by $e_i(\M + \X_1) - e_i(\M' + \X'_1)$ for each $i$. We expand this expression in the quotient of $\RS_{a, b}$ by the relations (1) and (4):

\begin{align*}
	e_i(\X_1 + \M) - e_i(\X'_1 + \M') & = \sum_{j + k = i} e_j(\X_1)e_k(\M) - e_j(\X'_1)e_k(\M') \\
	& = \sum_{j + k + \ell = i} (-1)^\ell e_j(\X_1)e_k(\X_2)h_\ell(\B) - (-1)^\ell e_j(\X'_1)e_k(\X'_2)h_\ell(\B) \\
	& = \sum_{j + k + \ell = i} (-1)^\ell h_\ell(\B) \left( e_j(\X_1) e_k(\X_2) - e_j(\X'_1) e_k(\X'_2) \right) \\
	& = \sum_{j + \ell = i} (-1)^\ell h_\ell(\B) \left( e_j(\X_1 + \X_2) - e_j(\X'_1 + \X'_2) \right).
\end{align*}

In particular, for each $1 \leq i \leq a+b-k$, we have

\[
\left(e_i(\X_1 + \M) - e_i(\X'_1 + \M') \right) - \left(e_i(\X_1 + \X_2) - e_i(\X'_1 + \X'_2) \right) \in \langle e_j(\X_1 + \X_2) - e_j(\X'_1 + \X'_2) \ | \ j \leq i \rangle.
\]
It follows by induction on $i$ that, after passing to the quotient by (1) and (4), we have

\[
\langle e_i(\X_1 + \M) - e_i(\X'_1 + \M') \ | \ 1 \leq i \leq a+b-k \rangle =  \langle e_i(\X_1 + \X_2) - e_i(\X'_1 + \X'_2) \ | \ 1 \leq i \leq a+b-k \rangle.
\]

\end{proof}

\subsection{Unreduced Colored HOMFLY Homology}

In this section we recall the construction of unreduced colored HOMFLY homology via singular Soergel bimodules from \cite{WW17, MSV11}.  We begin by describing the Rickard complex of singular Soergel bimodules associated to a colored braid diagram.

\begin{defn}
A \textit{colored braid diagram} $\beta$ is a braid diagram with colors associated to each strand. The \textit{bottom} and \textit{top} compositions of $\beta$ are the sequences of colors, read from left to right, assigned to the bottom and top of the diagram, respectively.
\end{defn}

\begin{ex}
The following is a colored braid diagram with bottom composition $(4,4,3)$ and top composition $(4,3,4)$:

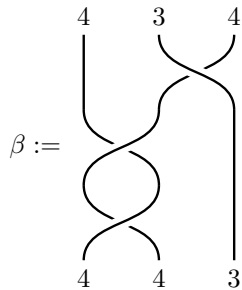
\begin{gather*}
\beta := 
\begin{tikzpicture}[anchorbase]
	\draw[webs] (1,0) to[out=90,in=270] (0,1);
	\draw[line width=5pt, color=white] (0,0) to[out=90,in=270] (1,1);
	\draw[webs] (0,0) to[out=90,in=270] (1,1);
	\draw[webs] (2,-1) node[below]{$3$} to (2,1);
	\draw[webs] (1,1) to[out=90,in=270] (2,2) node[above]{$4$};
	\draw[line width=5pt, color=white] (2,1) to[out=90,in=270] (1,2);
	\draw[webs] (2,1) to[out=90,in=270] (1,2) node[above]{$3$};	
	\draw[webs] (1,-1) node[below]{$4$} to[out=90,in=270] (0,0);
	\draw[line width=5pt,color=white] (0,-1) to[out=90,in=270] (1,0);
	\draw[webs] (0,-1) node[below]{$4$} to[out=90,in=270] (1,0);
	\draw[webs] (0,1) to (0,2) node[above]{$4$};
\end{tikzpicture}
\end{gather*}
\end{ex}

\begin{defn} \label{def: Rickard_complex_crossing}
    Fix $a, b \geq 0$, and set $m := \min(a, b)$. Then the \textit{2-strand positive Rickard complex} $C_{a, b}$ is the bounded complex of singular Soergel bimodules

	\begin{center}
	\begin{tikzcd}
	C_{a,b} := W_0 \arrow[r] & q^{-1}t W_1 \arrow[r] & \dots \arrow[r] & q^{-m+1}t^{m-1} W_{m-1} \arrow[r] & q^{-m}t^m W_m \in \CS(\SSBim_{a, b}^{b, a}).
	\end{tikzcd}
	\end{center}
	Similarly, the \textit{2-strand negative Rickard complex} $C_{a, b}^\vee$ is the bounded complex of singular Soergel bimodules
	
	\begin{center}
	\begin{tikzcd}
	C_{a,b}^{\vee} := q^mt^{-m} W_m \arrow[r] & q^{m-1}t^{-m+1} W_{m-1} \arrow[r] & \dots \arrow[r] & qt^{-1} W_1 \arrow[r] & W_0 \in \CS(\SSBim_{a, b}^{b, a}).
	\end{tikzcd}
	\end{center}
	
	We will not be concerned with the exact form of the differentials in this complex; see e.g. \cite{Wed19} for details.
\end{defn}

\begin{defn}
Let $\beta$ be a colored braid word. Assign to each colored crossing a $2$-strand positive or negative Rickard complex as follows:

\begin{align*}
F
\begin{tikzpicture}[anchorbase]
\draw[webs] (1,0) node[below]{$b$} to[out=90,in=270] (0,1);
\draw[line width=5pt,color=white] (0,0) to[out=90,in=270] (1,1);
\draw[webs] (0,0) node[below]{$a$} to[out=90,in=270] (1,1);
\draw (-.25,-.5) to[out=125,in=235] (-.25,1.25);
\draw (1.25,-.5) to[out=55,in=305] (1.25,1.25);
\end{tikzpicture}
= C_{a,b}; \quad \quad F
\begin{tikzpicture}[anchorbase]
\draw[webs] (0,0) node[below]{$a$} to[out=90,in=270] (1,1);
\draw[line width=5pt,color=white] (1,0) to[out=90,in=270] (0,1);
\draw[webs] (1,0) node[below]{$b$} to[out=90,in=270] (0,1);
\draw (-.25,-.5) to[out=125,in=235] (-.25,1.25);
\draw (1.25,-.5) to[out=55,in=305] (1.25,1.25);
\end{tikzpicture}
= C^\vee_{a,b}.
\end{align*}

The \textit{Rickard complex} $F(\beta)$ assigned to $\beta$ is the complex obtained from the tensor product of all $2$-strand Rickard complexes associated to crossings of $\beta$, arranged to respect the tensor products corresponding to horizontal and vertical concatenation in the diagrammatic description of $\SSBim$.
\end{defn}

The chain bimodules appearing in a general Rickard complex will be tensor products of ladder webs $W_k$ assigned to each crossing. We give some terminology for such bimodules next.

\begin{defn} \label{def: resolution}
Given a colored braid diagram $\beta$, fix some ordering $\{c_1, \dots, c_n\}$ of the crossings of $\beta$. Denote the two colors of the strands entering $c_i$ by $a_i, b_i$. A \textit{resolution} of $\beta$ is a choice $\vec{k} = (k_1, \dots, k_n)$ of a non-negative integer $0 \leq k_i \leq \min(a_i, b_i)$ for each $i$. For each resolution $\vec{k}$ of $\beta$, we denote by $W_{\vec{k}}$ the corresponding tensor product of ladder webs which appears in the Rickard complex $F(\beta)$.
\end{defn}

Suppose $p$ is a point on a colored braid diagram $\beta$. Denote the color of the strand containing $p$ by $k$. For each resolution $\vec{k}$ of $\beta$, the corresponding bimodule $W_{\vec{k}}$ will contain a tensor factor of $\Sym(\X_p)$, where $\X_p$ is an alphabet of size $k$. This describes an action of $\Sym(\X_p)$ on the Rickard complex $F(\beta)$ given by multiplication in this tensor factor. Many of these actions turn out to be homotopic.

\begin{lem} \label{lem: dot-sliding}
The actions of $e_i(\X_1) - e_i(\X'_2)$ and $e_i(\X_2) - e_i(\X'_1)$ on a $2$-strand Rickard complex are nullhomotopic. Moreover, these homotopies can be chosen to pairwise anti-commute.
\end{lem}

\begin{proof}
Standard; see e.g. \cite{HRW24}, Lemma 4.20.
\end{proof}

We refer to the homotopies of Lemma \ref{lem: dot-sliding} as \textit{dot-sliding} homotopies. By taking linear combinations of dot-sliding homotopies, we see that the actions of $f(\X_p)$ and $f(\X_q)$ are homotopic for any two points $p, q$ on a $k$-labeled strand of $\beta$ and any symmetric polynomial $f$ in $k$ variables.

\begin{notn}
Let $\aa$ and $\bb$ be the top and bottom compositions of $\beta$. Then in particular, there are actions of $\Sym^\aa(\X)$ and $\Sym^\bb(\X)$ on $F(\beta)$; this is exactly the usual bimodule action in $\SSBim_\aa^\bb$. We label the associated alphabets on the top and bottom by $\X^{ext}_i$ and $\X'^{ext}_i$, respectively, where the index $i$ increases from left to right. In this convention, we have $|\X^{ext}_i| = a_i$ and $|\X'^{ext}_i| = b_i$.
\end{notn}

There is a natural inclusion of $\SSBim_\aa^\bb$ into $D^b(\SSBim_\aa^\bb)$ given by considering a bimodule as a complex concentrated in Hochschild degree $0$. This extends to an inclusion on the corresponding homotopy categories, so we can naturally consider the Rickard complex $F(\beta)$ as living in $K^b(D^b(\SSBim_\aa^\bb))$.

\begin{notn}
Given two alphabets $\X_i$, $\X_j$ of the same size $k$, we abbreviate the Koszul (in the Hochschild or homological directions) complex on the difference of corresponding elementary symmetric polynomials in these two alphabets as

\[
K(C; \X_i - \X_j) := K(C; \{e_r(\X_i) - e_r(\X_j) \ | \ 1 \leq r \leq k \}).
\]
We similarly write $K(C; \X_i)$ for the Koszul complex on the elementary symmetric polynomials in an alphabet $\X_i$ and extend this notation to multiple alphabets and multiple differences of alphabets in the natural way.
\end{notn}

\begin{defn} \label{def: hoch_hom}
Given a composition $\aa = (a_1, \dots, a_n) \vdash N$ and a complex $C \in K^b(D^b(\SSBim_\aa^\aa))$, the \textit{Hochschild homology} of $C$ is the Koszul complex

\[
HH(C) := K_a(C; \{\X^{ext}_i - \X'^{ext}_i | 1 \leq i \leq n\}).
\]

We denote by $H(HH(C))$ the $\Z_a \times \Z_q \times \Z_t$-graded $R$-module obtained from $HH(C)$ by first taking homology in the Hochschild direction to obtain a complex of $\Z_a \times \Z_q$-graded $R$-modules, then taking homology of the result.
\end{defn}

We can consider $H(HH(F(\beta)))$ for any colored braid diagram $\beta$ with equal top and bottom compositions. Much of the module structure here degenerates. Since the actions of symmetric polynomials on any given braid strand are homotopic, these actions are identified on the nose after taking homology in the $t$-direction. Moreover, taking homology in the $a$-direction of the Koszul complex from Definition \ref{def: hoch_hom} identifies on the nose the actions of $\Sym(\X^{ext}_i)$ and $\Sym(\X'^{ext}_i)$ for each $i$. In total, we're left with an action of the following ring.

\begin{defn} \label{def: link_symms}
Given a colored braid diagram $\beta$ with equal top and bottom composition $\aa = (a_1, \dots, a_n)$, denote by $\pi_0(\hat{\beta})$ the set of orbits of $\{1, \dots, n\}$ under the permutation underlying $\beta$. We set

\[
\Sym(\X_{\hat{\beta}}) := \bigotimes_{[i] \in \pi_0(\hat{\beta})} \Sym(\X^{ext}_i).
\]
\end{defn}

The following theorem is due to Khovanov \cite{Kh07} and Khovanov--Rozansky \cite{KhR08} in the uncolored setting and Mackaay-Stošić-Vaz \cite{MSV11} and Webster--Williamson \cite{WW17} in the colored setting.

\begin{thm} \label{thm: hhh_invariance}
Let $\mathcal{L}$ be a colored link presented as the closure of a colored braid diagram $\beta$. Then $HHH(\mathcal{L}) := H(HH(F(\beta)))$ is independent of the choice of $\beta$ up to $\Sym(\X_{\hat{\beta}})$-equivariant isomorphism. We call $HHH(\L)$ the \textit{(unreduced) colored HOMFLY homology} of $\mathcal{L}$.
\end{thm}

\subsection{Reduced Rickard Complexes}
We now turn to giving a construction of reduced colored HOMFLY homology which generalizes the results of \cite{Wed19} to arbitrary colored links reduced at an arbitrary component. We mostly follow the notation from \cite{Wed19}.

To begin, we fix a colored link $\L$, a colored braid diagram $\beta$ with $\hat{\beta} = \L$, and a marked point $p \in \beta$ away from the crossings. We write $\aa = (a_1, \dots, a_n)$ for the composition labeling the external edges of $\beta$; since $\beta$ closes to a link, this composition is the same at the top and bottom. We write $i \sim j$ if the $i^{th}$ and $j^{th}$ braid strand close to the same link component.

Reading the strands from left to right along the bottom of $\beta$, we denote the index of the strand containing the \textbf{m}arked point $p$ by $m$. We label the alphabet associated to the segment containing $p$ as $\X_p$, so that $|\X_p| = a_m$.

\begin{ex} \label{ex: reduction_data}
We depict an example of the required data below:

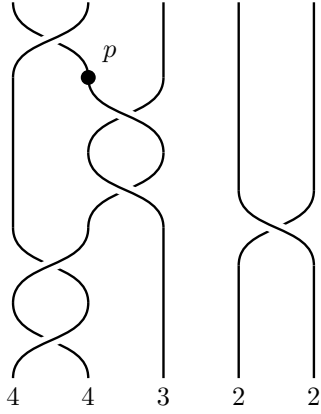
\begin{center}
\begin{tikzpicture}[anchorbase]
	\draw[webs] (1,0) to[out=90,in=270] (0,1);
	\draw[line width=5pt, color=white] (0,0) to[out=90,in=270] (1,1);
	\draw[webs] (0,0) to[out=90,in=270] (1,1);
	\draw[webs] (2,-1) node[below]{$3$} to (2,1);
	\draw[webs] (1,1) to[out=90,in=270] (2,2);
	\draw[line width=5pt, color=white] (2,1) to[out=90,in=270] (1,2);
	\draw[webs] (2,1) to[out=90,in=270] (1,2);
	\draw[webs] (1,2) to[out=90,in=270] (2,3);
	\draw[line width=5pt, color=white] (2,2) to[out=90,in=270] (1,3);
	\draw[webs] (2,2) to[out=90,in=270] (1,3);
	\draw[webs] (1,3) to[out=90,in=270] (0,4);
	\draw[line width=5pt, color=white] (0,3) to[out=90,in=270] (1,4);
	\draw[webs] (0,3) to[out=90,in=270] (1,4);
	\draw[webs] (3,-1) node[below]{$2$} to (3, .5);
	\draw[webs] (4,-1) node[below]{$2$} to (4,.5);
	\draw[webs] (3,.5) to[out=90,in=270] (4,1.5);
	\draw[color=white, line width=5pt] (4,.5) to[out=90,in=270] (3,1.5);
	\draw[webs] (4,.5) to[out=90,in=270] (3,1.5);
	\draw[webs] (3,1.5) to (3,4);
	\draw[webs] (4,1.5) to (4,4);
	\draw[webs] (0,1) to (0,3);
	\draw[webs] (2,3) to (2,4);
	
	\draw[webs] (1,-1) node[below]{$4$} to[out=90,in=270] (0,0);
	\draw[line width=5pt,color=white] (0,-1) to[out=90,in=270] (1,0);
	\draw[webs] (0,-1) node[below]{$4$} to[out=90,in=270] (1,0);
	 
	\node[circle,fill=black,inner sep=2pt,label={north east:$p$}] at (1,3) {};;
\end{tikzpicture}
\end{center}
Here $\aa = (4, 4, 3, 2, 2) \vdash 15$, $m = 1$, and $|\X_p| = a_1 = 4$. Notice that $\beta$ splits as a disjoint union of two braids and that $a_m$ is not minimal among the strand labels of $\beta$.
\end{ex}

\begin{defn}
Recall the notation $\RS_{a, b}$, $\IS_{a, b}$, and $\JS_k$ from Definition \ref{def: new_web_rels}. For each $(a,b)$-crossing, set  $\VS_{a,b} := \RS_{a,b}/\IS_{a,b}$ and denote by $\Scal_{a,b}$ the subring of $\VS_{a,b}$ generated by $\Sym(\X_1 | \X_2 | \X'_1 | \X'_2)$.
\end{defn}

\begin{lem} \label{lem: vs_poly_alg}
$\Scal_{a,b}$ is isomorphic to a polynomial algebra over $R$, and $\VS_{a, b}$ is isomorphic to a polynomial algebra over $\Scal_{a, b}$.
\end{lem}

\begin{proof}
Using Corollary \ref{cor: solving_for_alphabets}, we can rewrite the first, second, and third family of relations in $\IS_{a,b}$ to express the generators $e_i(\M)$, $e_i(\M')$, and $e_i(\F)$ of $\RS_{a,b}$, respectively, in terms of the alphabets $\X_1, \X_2, \X_2'$, and $\B$ as follows:

\begin{align*}
e_i(\M) & = e_i(\X_2 - \B); \\
e_i(\M') & = e_i(\X'_2 - \B); \\
e_i(F) & = e_i(\X_1 + \M).
\end{align*}
Making these substitutions, we obtain an isomorphism

\[
\VS_{a,b} \cong \Sym(\X_1 | \X_2 | \X'_1 | \X'_2 | \B) / \langle e_i(\X_1 + \X_2) - e_i(\X'_1 + \X'_2) \ | \ 1 \leq i \leq \max(a,b) \rangle.
\]

Now if $b \geq a$, we can again apply Corollary \ref{cor: solving_for_alphabets} to the remaining relations to express either the generators $e_i(\X'_2)$ in terms of the alphabets $\X_1, \X_2$, and $\X'_1$ or the generators $e_i(\X_1)$ in terms of the alphabets $\X_2, \X'_1$, and $\X'_2$, giving isomorphisms

\begin{equation} \label{eq: scal_poly}
\VS_{a,b} \cong \Sym(\X_1 | \X_2 | \X'_1 | \B) \cong \Sym(\X_2 | \X_1' | \X'_2 | \B).
\end{equation}
Then $\Scal_{a,b} \cong \Sym(\X_1 | \X_2 | \X'_1) \cong \Sym(\X_2 | \X'_1 | \X'_2)$, so

\[
\VS_{a,b} \cong \Scal_{a,b}[e_i(\B) \ | \ 1 \leq i \leq b].
\]
Otherwise $a > b$, and one can write either $e_i(\X_2)$ or $e_i(\X'_1)$ in terms of the other alphabets; then $\VS_{a,b}$ is still a polynomial ring over $\Scal_{a,b}$ in the alphabet $\{e_i(\B)\}$.
\end{proof}

\begin{prop} \label{prop: ideal_regular_crossing}
The generators of $\JS_k$ form a regular sequence in $\VS_{a,b}$.
\end{prop}

\begin{proof}
Essentially identical to the proof of Proposition \ref{prop: web_quot_regular}.
\end{proof}

\begin{cor}
There is a quasi-isomorphism

\[
W_k \cong K(\VS_{a, b}; \JS_k) \in D^b(\Scal_{a,b}-\Mod{}).
\]
\end{cor}

We now extend the local notions defined above to the entire braid diagram $\beta$:

\begin{defn}
Let $\RS_\beta$ be the tensor product of the rings $\RS_{a,b}$ associated to every crossing, in which the alphabets $\X_i$ and $\X'_i$ which are assigned to the same segment between crossings of $\beta$ are identified. Similarly let $\IS_\beta \subset \RS_\beta$ be the product of all the ideals $\IS_{a,b}$ assigned to crossings, and set $\VS_\beta := \RS_\beta/\IS_\beta$. Denote by $\Scal_\beta$ the subring of $\VS_\beta$ generated by (symmetric polynomials in) the alphabets $\X_i, \X'_i$ assigned to each crossing. For each resolution $\vec{k}$ of $\beta$, let $\JS_{\vec{k}}$ denote the product of all the corresponding ideals $\JS_{k_i}$.
\end{defn}

The facts presented above for the local data have immediate analogues for the entire braid. In formulating a global version of Lemma \ref{lem: vs_poly_alg}, it will be useful to take some care to avoid eliminating certain alphabets.

\begin{prop} \label{prop: scal_poly_beta}
$\Scal_\beta$ is isomorphic to a polynomial algebra over $\Sym(\X'^{ext}_1 | \dots | \X'^{ext}_{m - 1} | \X_p | \X'^{ext}_{m + 1} | \dots | \X'^{ext}_n)$, and $\VS_\beta$ is isomorphic to a polynomial algebra over $\Scal_\beta$.
\end{prop}

\begin{proof}
To see that $\Scal_\beta$ is a polynomial algebra over the appropriate alphabets, we can run the argument of Proposition \ref{lem: vs_poly_alg} at each crossing; the only potential obstructions are that we must eliminate a different alphabet at each crossing and that we can never eliminate the alphabet $\X_p$ or any of the alphabets $\X'^{ext}_i$ for $i \neq m$.

To make such a choice of alphabets to eliminate, order the strands of $\beta$ by decreasing color. At each stage, if the next strand does not contain $p$, travel upwards from the bottom of that strand, selecting at each crossing which has not yet been treated the alphabet on that strand which enters that crossing from above. If the next strand does contain $p$, then instead travel outwards from $p$, selecting at each crossing which has not yet been treated the alphabet on that strand farthest from $p$. This produces an alphabet associated to each crossing without repetition, ensures that $\X_p$ and $\X'^{ext}_i$ for $i \neq m$ are never chosen, and ensures that the alphabet associated to each crossing lies on a strand of greatest color passing through that crossing.

Eliminating these alphabets using the $4^{th}$ family of relations in $\IS_\beta$ presents $\VS_\beta$ as a polynomial algebra generated by $\Sym(\X_p)$, $\Sym(\X'^{ext}_i)$ for $i \neq m$, potentially other $\Sym(\X_i)$, and $\Sym(\B)$ at each crossing. We obtain an analogous description of $\Scal_\beta$ by forgetting the generators of $\Sym(\B)$ at each crossing.
\end{proof}

\begin{ex}
In the context of Example \ref{ex: reduction_data}, order the strands from left to right along the bottom of $\beta$. Then the above algorithm eliminates the following alphabets (marked with {\color{red}{red}} dots):

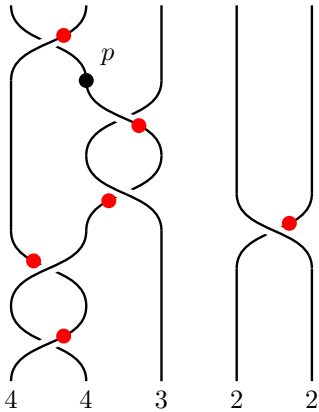
\begin{center}
\begin{tikzpicture}[anchorbase]
	\draw[webs] (1,0) to[out=90,in=270] (0,1);
	\draw[line width=5pt, color=white] (0,0) to[out=90,in=270] (1,1);
	\draw[webs] (0,0) to[out=90,in=270] (1,1);
	\draw[webs] (2,-1) node[below]{$3$} to (2,1);
	\draw[webs] (1,1) to[out=90,in=270] (2,2);
	\draw[line width=5pt, color=white] (2,1) to[out=90,in=270] (1,2);
	\draw[webs] (2,1) to[out=90,in=270] (1,2);
	\draw[webs] (1,2) to[out=90,in=270] (2,3);
	\draw[line width=5pt, color=white] (2,2) to[out=90,in=270] (1,3);
	\draw[webs] (2,2) to[out=90,in=270] (1,3);
	\draw[webs] (1,3) to[out=90,in=270] (0,4);
	\draw[line width=5pt, color=white] (0,3) to[out=90,in=270] (1,4);
	\draw[webs] (0,3) to[out=90,in=270] (1,4);
	\draw[webs] (3,-1) node[below]{$2$} to (3, .5);
	\draw[webs] (4,-1) node[below]{$2$} to (4,.5);
	\draw[webs] (3,.5) to[out=90,in=270] (4,1.5);
	\draw[color=white, line width=5pt] (4,.5) to[out=90,in=270] (3,1.5);
	\draw[webs] (4,.5) to[out=90,in=270] (3,1.5);
	\draw[webs] (3,1.5) to (3,4);
	\draw[webs] (4,1.5) to (4,4);
	\draw[webs] (0,1) to (0,3);
	\draw[webs] (2,3) to (2,4);
	
	\draw[webs] (1,-1) node[below]{$4$} to[out=90,in=270] (0,0);
	\draw[line width=5pt,color=white] (0,-1) to[out=90,in=270] (1,0);
	\draw[webs] (0,-1) node[below]{$4$} to[out=90,in=270] (1,0);
	 
	\node[circle,fill=black,inner sep=2pt,label={north east:$p$}] at (1,3) {};;
	
	\node[circle,fill=red,inner sep=2pt] at (.7,-.4) {};
	\node[circle,fill=red,inner sep=2pt] at (.3,.6) {};
	\node[circle,fill=red,inner sep=2pt] at (1.3,1.4) {};
	\node[circle,fill=red,inner sep=2pt] at (1.7,2.4) {};
	\node[circle,fill=red,inner sep=2pt] at (.7,3.6) {};
	\node[circle,fill=red,inner sep=2pt] at (3.7,1.1) {};
\end{tikzpicture}
\end{center}
Observe that an alphabet is eliminated at every crossing, that alphabet always lies on a strand of maximal color entering the crossing, alphabets are not repeated, $\X_p$ is not chosen, and $\X'^{ext}_i$ is never chosen for $i \neq m$ (nor in this example when $i = m$, though this is not required).

\end{ex}

\begin{prop} \label{prop: ideal_regular_braid}
$W_{\vec{k}} \cong \VS_\beta / \JS_{\vec{k}}$, and the generators of $\JS_{\vec{k}}$ form a regular sequence in $\VS_\beta$.
\end{prop}

\begin{proof}
Completely analogous to the local case (Proposition \ref{prop: ideal_regular_crossing}).
\end{proof}

\begin{cor}
There is a quasi-isomorphism 

\[
W_{\vec{k}} \cong K(\VS_\beta; \JS_{\vec{k}}) \in D^b(\Scal_\beta-\Mod{}).
\]
\end{cor}

\begin{defn}
Set $\overline{\Scal_\beta(p)} :=\Scal_\beta/\Sym(\X_p)$.
\end{defn}

Because the differential in the Rickard complex associated to each crossing consists of bimodule homomorphisms, the differential on the complex $F(\beta)$ associated to the whole braid is $\Scal_\beta$-linear. We can therefore consider $F(\beta)$ as a complex in $K^b(D^b(\Scal_\beta-\Mod{}))$ in the following definition.

\begin{defn}
The \textit{reduced Rickard complex} of $\beta$ at $p$ is the complex 

\[
F_p(\beta) := F(\beta) \tensor[\Scal_\beta, der] \overline{\Scal_\beta(p)} \in K^b(D^b(\Scal_{\beta}-\Mod{})).
\]

\end{defn}

We close this section by establishing an invariance result for the reduced Rickard complex; see Theorem \ref{thm: red_rickard_invariance} below.

\begin{prop} \label{prop: quot_for_cone_rickard}
$F_p(\beta) \simeq K_t(F(\beta); \X_p) \in K^b(D^b(\overline{\Scal_\beta(p)}-\Mod{}))$.
\end{prop}

\begin{proof}
Follows from Proposition \ref{prop: scal_poly_beta} and repeated applications of Proposition \ref{prop: quot_for_cone}.
\end{proof}

\begin{defn}
For each strand $\lambda$ of $\beta$, denote by $\overline{\Scal_\beta(\lambda)}$ the subring of $\Scal_\beta$ generated by (symmetric polynomials in) alphabets which do not lie on the strand $\lambda$.
\end{defn}

By Proposition \ref{prop: scal_poly_beta}, we can view $\overline{\Scal_\beta(p)}$ as the subring of $\Scal_\beta$ generated by (symmetric polynomials in) all alphabets except $\X_p$. It follows that there's a natural inclusion $\overline{\Scal_\beta(\lambda)} \hookrightarrow \overline{\Scal_\beta(p)}$ for each $p \in \lambda$. Given two points $p, p' \in \lambda$, we have a commutative diagram of inclusions

\begin{center}
\begin{tikzcd}[sep = large]
& \Scal_\beta & \\
\overline{\Scal_\beta(p)} \arrow[hookrightarrow]{ur} & & \overline{\Scal_\beta(p')} \arrow[hookrightarrow]{ul} \\
& \overline{\Scal_\beta(\lambda)} \arrow[hookrightarrow]{ul} \arrow[hookrightarrow]{ur} &
\end{tikzcd}
\end{center}
This in turn furnishes a commutative diagram of functors given by restriction of scalars along these inclusions, and we implicitly apply these restrictions in what follows.

\begin{thm} \label{thm: red_rickard_invariance}
Given a braid strand $\lambda$ of $\beta$ and two points $p, p' \in \lambda$, there is a homotopy equivalence of reduced Rickard complexes

\[
F_p(\beta) \simeq F_{p'}(\beta) \in K^b(D^b(\overline{\Scal_\beta(\lambda)}-\Mod{})).
\]
\end{thm}

\begin{proof}
By Proposition \ref{prop: quot_for_cone_rickard}, we have homotopy equivalences $F_p(\beta) \simeq K_t(F(\beta); \X_p)$ and $F_{p'}(\beta) \simeq K_t(F(\beta); \X_{p'})$ in $K^b(D^b(\overline{\Scal_\beta(\lambda)}))$. Since $p$ and $p'$ lie on the same strand, the endomorphisms $e_i(\X_p)$ and $e_i(\X_{p'})$ of $F(\beta)$ are homotopic for each $i$ via a graded commutative family of homotopies. The desired result then follows from repeated application of Lemma \ref{lem: kosz_homotopic}.
\end{proof}

\begin{rem}
We wish to stress that the homotopy equivalence of Theorem \ref{thm: red_rickard_invariance} occurs in $K^b(D^b(\overline{\Scal_\beta(\lambda)}-\Mod{}))$ and \textit{not} in $K^b(D^b(\Scal_\beta-\Mod{}))$. While each reduced Rickard complex $F_p(\beta)$ and $F_{p'}(\beta)$ most naturally lives in $K^b(D^b(\Scal_\beta-\Mod{}))$, the equivalence between them does not respect the full action of $\Scal_\beta$. For example, $\Sym(\X_p)$ acts by $0$ on each chain module on the left while it acts nontrivially on each chain module on the right.

This presents a choice: does one define the reduced Rickard complex as reduced at a \textit{point} with an action of $\Scal_\beta$ or as reduced at a \textit{strand} with an action of $\overline{\Scal_\beta(\lambda)}$? We take the former approach: we are ultimately interested not in braid invariants but in \textit{link} invariants, and the action of the larger ring $\Scal_\beta$ is essential in extracting Hochschild homology below. We would be very interested to see a presentation of this invariant in which Rickard complexes are reduced at a strand before taking Hochschild homology.
\end{rem}

\subsection{Reduced Column-Colored Homology} \label{subsec: reduced_homology}

\begin{defn}
We denote by $HH_{\neg p}$ the functor on $D^b(\SSBim_{\aa}^{\aa})$ given by

\[
HH_{\neg p}(C) := K_a\left(C; \{\X^{ext}_i - \X'^{ext}_i \ | i \neq m\} \right).
\]
$HH_{\neg p}$ is additive, so it extends automatically to $K^b(D^b(\SSBim_{\aa}^{\aa}))$.
\end{defn}

Graphically, $HH_{\neg p}(F_p(\beta))$ corresponds to closing the braid $\beta$ at all strands \textit{except} the strand $\lambda$ containing $p$, which we leave open at the bottom. In our running example, the reader should have the following picture in mind:

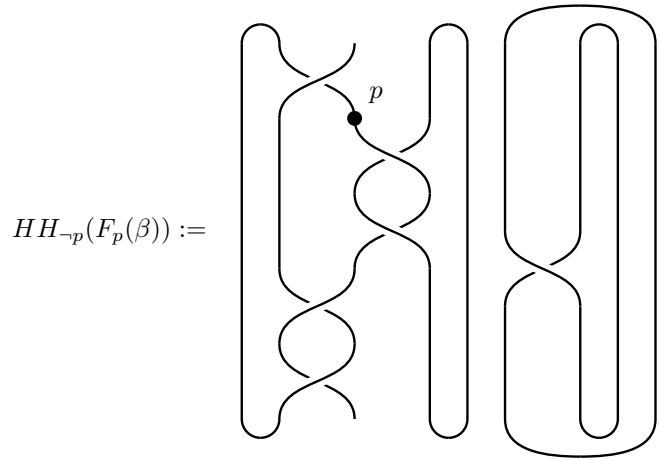
\begin{gather*}
HH_{\neg p}(F_p(\beta)) := \quad
\begin{tikzpicture}[anchorbase]
	\draw[webs] (1,0) to[out=90,in=270] (0,1);
	\draw[line width=5pt, color=white] (0,0) to[out=90,in=270] (1,1);
	\draw[webs] (0,0) to[out=90,in=270] (1,1);
	\draw[webs] (2,-1) to (2,1);
	\draw[webs] (1,1) to[out=90,in=270] (2,2);
	\draw[line width=5pt, color=white] (2,1) to[out=90,in=270] (1,2);
	\draw[webs] (2,1) to[out=90,in=270] (1,2);
	\draw[webs] (1,2) to[out=90,in=270] (2,3);
	\draw[line width=5pt, color=white] (2,2) to[out=90,in=270] (1,3);
	\draw[webs] (2,2) to[out=90,in=270] (1,3);
	\draw[webs] (1,3) to[out=90,in=270] (0,4);
	\draw[line width=5pt, color=white] (0,3) to[out=90,in=270] (1,4);
	\draw[webs] (0,3) to[out=90,in=270] (1,4);
	\draw[webs] (3,-1) to (3, .5);
	\draw[webs] (4,-1) to (4,.5);
	\draw[webs] (3,.5) to[out=90,in=270] (4,1.5);
	\draw[color=white, line width=5pt] (4,.5) to[out=90,in=270] (3,1.5);
	\draw[webs] (4,.5) to[out=90,in=270] (3,1.5);
	\draw[webs] (3,1.5) to (3,4);
	\draw[webs] (4,1.5) to (4,4);
	\draw[webs] (0,1) to (0,3);
	\draw[webs] (2,3) to (2,4);	
	\draw[webs] (1,-1) to[out=90,in=270] (0,0);
	\draw[line width=5pt,color=white] (0,-1) to[out=90,in=270] (1,0);
	\draw[webs] (0,-1) to[out=90,in=270] (1,0);
	\node[circle,fill=black,inner sep=2pt,label={north east:$p$}] at (1,3) {};;
	\draw[webs] (4,-1) to[out=270,in=180] (4.25,-1.25);
	\draw[webs] (4.25,-1.25) to[out=0,in=270] (4.5,-1);
	\draw[webs] (4.5,-1) to (4.5,4);
	\draw[webs] (4.5,4) to[out=90,in=0] (4.25,4.25);
	\draw[webs] (4.25,4.25) to[out=180,in=90] (4,4);
	\draw[webs] (3,-1) to[out=270,in=180] (4,-1.5);
	\draw[webs] (4,-1.5) to[out=0,in=270] (5,-1);
	\draw[webs] (5,-1) to (5,4);
	\draw[webs] (5,4) to[out=90,in=0] (4,4.5);
	\draw[webs] (4,4.5) to[out=180,in=90] (3,4);
	\draw[webs] (2,-1) to[out=270,in=180] (2.25,-1.25);
	\draw[webs] (2.25,-1.25) to[out=0,in=270] (2.5,-1);
	\draw[webs] (2.5,-1) to (2.5,4);
	\draw[webs] (2.5,4) to[out=90,in=0] (2.25,4.25);
	\draw[webs] (2.25,4.25) to[out=180,in=90] (2,4);
	\draw[webs] (0,-1) to[out=270,in=0] (-.25,-1.25);
	\draw[webs] (-.25,-1.25) to[out=180,in=270] (-.5,-1);
	\draw[webs] (-.5,-1) to (-.5,4);
	\draw[webs] (-.5,4) to[out=90,in=180] (-.25,4.25);
	\draw[webs] (-.25,4.25) to[out=0,in=90] (0,4);
\end{tikzpicture}
\end{gather*}

\begin{defn}
Given a component $\L_p$ of $\L$, we denote by $\overline{\Scal_\beta(\L_p)}$ the subring of $\Scal_\beta$ generated by (symmetric polynomials in) alphabets which do not lie on $\L_p$.
\end{defn}

\begin{prop} \label{prop: kosz_a_for_kosz_t_Rickard}
There is a homotopy equivalence in $K^b(D^b(\overline{\Scal_\beta(\L_p)}-\Mod{}))$

\[
HH_{\neg p}(F_p(\beta)) \simeq K_t\left( K_a \left( F(\beta); \{\X_i^{ext} - \X'^{ext}_i \}_{i \not \sim m} \right); \X_p \sqcup \{\X_i^{ext} - \X'^{ext}_i\}_{m \neq i \sim m} \right).
\]
\end{prop}

In words, Proposition \ref{prop: kosz_a_for_kosz_t_Rickard} states that among the Koszul complexes in the Hochschild direction which define $HH_{\neg p}$, those which lie on the same component of $\L$ as $p$ can be exchanged for Koszul complexes in the homological direction.

\begin{proof}
By Proposition \ref{prop: scal_poly_beta}, $\Scal_\beta$ is isomorphic to a polynomial algebra which contains among its generators $e_j(\X_p)$ and $e_j(\X'^{ext}_i)$ for $m \neq i \sim m$ and $1 \leq j \leq a_m$. By a change of variables, we can realize $\Scal_\beta$ as a polynomial algebra over the collection $\{e_j(\X^{ext}_i) - e_j(\X'^{ext}_i) \ | \ i \neq m, i \sim m, 1 \leq j \leq a_m\}$. The result then follows immediately from Corollary \ref{cor: cone_a_for_cone_t}.
\end{proof}

\begin{thm} \label{thm: point_invariance}
Suppose $p, p'$ lie on the same component $\L_p$ of $\L$. Then there is a homotopy equivalence

\[
HH_{\neg p}(F_p(\beta)) \simeq HH_{\neg p'}(F_{p'}(\beta)) \in K^b(D^b(\overline{\Scal_\beta(\L_p)}-\Mod{})).
\]
\end{thm}

\begin{proof}

Recall that we denote by $m$ the index of the strand of $\beta$ containing $p$; we will similarly denote by $m'$ the index of the strand containing $p'$. Since $p$ and $p'$ lie on the same component of $\mathcal{L}$, we have $m \sim m'$. Applying Proposition \ref{prop: kosz_a_for_kosz_t_Rickard} to both reduced Rickard complexes, we have

\begin{align*}
HH_{\neg p}(F_p(\beta)) & \simeq K_t\left( K_a \left( F(\beta); \{\X_i^{ext} - \X'^{ext}_i \}_{i \not \sim m} \right); \X_p \sqcup \{\X_i^{ext} - \X'^{ext}_i\}_{m \neq i \sim m} \right); \\
HH_{\neg p'}(F_{p'}(\beta)) & \simeq K_t\left( K_a \left( F(\beta); \{\X_i^{ext} - \X'^{ext}_i \}_{i \not \sim m'} \right); \X_{p'} \sqcup \{\X_i^{ext} - \X'^{ext}_i\}_{m \neq i \sim m'} \right).
\end{align*}

Since $m \sim m'$, the inner (Hochschild) Koszul complexes in these two expressions exactly agree. In the outer (homological) Koszul complexes, passage from $HH_{\neg p}(F_p(\beta))$ to $HH_{\neg p'}(F_{p'}(\beta))$ requires replacing $\X_p$ with $\X_{p'}$ and $\X^{ext}_{m'} - \X^{ext}_{m'}$ with $\X^{ext}_m - \X'^{ext}_m$. To make the former replacement, one can replace $\X_p$ with $\X_q$ where $q$ is any other point lying on the strand of $\beta$ containing $p$ using dot slides. If $p$ and $p'$ lie on the same strand, this suffices. Else, choose $q$ to be at the top of the strand. Then $\X_q = \X^{ext}_i$ for some $m \neq i \sim m$. By Corollary \ref{cor: row_ops}, we can replace $\X_q$ by $\X'^{ext}_i$ using the component $\X^{ext}_i - \X'^{ext}_i$ in the homological Koszul complex for $HH_{\neg p}(F_p(\beta))$. This alphabet lies on a different strand of $\beta$, and one can continue this procedure.

A similar procedure allows the replacement of $\X^{ext}_{m'} - \X'^{ext}_{m'}$ with $\pm(\X^{ext}_m - \X'^{ext}_m)$ by following a sequence of dot slides and row operations involving the differences $\X^{ext}_i - \X'^{ext}_i$. These row operations must never use the term $\X^{ext}_{m'} - \X'^{ext}_{m'}$, but this is easy to ensure by traveling ``outwards" from the segment responsible for this term in the closure along $\L_p$. Since $\Cone(f) = \Cone(-f)$ for any $f$, we are finished.
\end{proof}

The sequence of replacements in the proof of Theorem \ref{thm: point_invariance} is best illustrated graphically.

\begin{ex}
Consider the braid $\beta$ below with the indicated choices of marked points $p$ and $p'$:

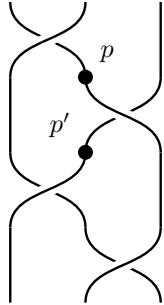
\begin{center}
\begin{tikzpicture}[anchorbase]
	\draw[webs] (2,0) to[out=90,in=270] (1,1);
	\draw[line width=5pt,color=white] (1,0) to[out=90,in=270] (2,1);
	\draw[webs] (1,0) to[out=90,in=270] (2,1);
	\draw[webs] (0,0) to (0,1);
	
	\draw[webs] (1,1) to[out=90,in=270] (0,2);
	\draw[line width=5pt,color=white] (0,1) to[out=90,in=270] (1,2);
	\draw[webs] (0,1) to[out=90,in=270] (1,2);
	\draw[webs] (2,1) to (2,2);
	
	\draw[webs] (1,2) to[out=90,in=270] (2,3);
	\draw[cover] (2,2) to[out=90,in=270] (1,3);
	\draw[webs] (2,2) to[out=90,in=270] (1,3);
	\draw[webs] (0,2) to (0,3);
	
	\draw[webs] (1,3) to[out=90,in=270] (0,4);
	\draw[cover] (0,3) to[out=90,in=270] (1,4);
	\draw[webs] (0,3) to[out=90,in=270] (1,4);
	\draw[webs] (2,3) to (2,4);
	
	\node[circle,fill=black,inner sep=2pt,label={north east:$p$}] at (1,3) {};;
	\node[circle,fill=black,inner sep=2pt,label={north west:$p'$}] at (1,2) {};;
\end{tikzpicture}
\end{center}
Then $HH_{\neg p}(F_p(\beta))$ is the following closure:

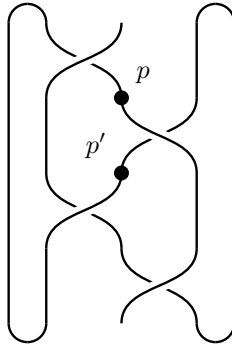
\begin{center}
\begin{tikzpicture}[anchorbase]
	\draw[webs] (2,0) to[out=90,in=270] (1,1);
	\draw[line width=5pt,color=white] (1,0) to[out=90,in=270] (2,1);
	\draw[webs] (1,0) to[out=90,in=270] (2,1);
	\draw[webs] (0,0) to (0,1);
	
	\draw[webs] (1,1) to[out=90,in=270] (0,2);
	\draw[line width=5pt,color=white] (0,1) to[out=90,in=270] (1,2);
	\draw[webs] (0,1) to[out=90,in=270] (1,2);
	\draw[webs] (2,1) to (2,2);
	
	\draw[webs] (1,2) to[out=90,in=270] (2,3);
	\draw[cover] (2,2) to[out=90,in=270] (1,3);
	\draw[webs] (2,2) to[out=90,in=270] (1,3);
	\draw[webs] (0,2) to (0,3);
	
	\draw[webs] (1,3) to[out=90,in=270] (0,4);
	\draw[cover] (0,3) to[out=90,in=270] (1,4);
	\draw[webs] (0,3) to[out=90,in=270] (1,4);
	\draw[webs] (2,3) to (2,4);
	
	\node[circle,fill=black,inner sep=2pt,label={north east:$p$}] at (1,3) {};;
	\node[circle,fill=black,inner sep=2pt,label={north west:$p'$}] at (1,2) {};;
	
	\draw[webs] (2,4) to[out=90,in=180] (2.25,4.25);
	\draw[webs] (2.25,4.25) to[out=0,in=90] (2.5,4);;
	\draw[webs] (2.5,4) to (2.5,0);
	\draw[webs] (2.5,0) to[out=270,in=0] (2.25,-.25);
	\draw[webs] (2.25,-.25) to[out=180,in=270] (2,0);
	
	\draw[webs] (0,4) to[out=90,in=0] (-.25,4.25);
	\draw[webs] (-.25,4.25) to[out=180,in=90] (-.5,4);;
	\draw[webs] (-.5,4) to (-.5,0);
	\draw[webs] (-.5,0) to[out=270,in=180] (-.25,-.25);
	\draw[webs] (-.25,-.25) to[out=0,in=270] (0,0);
\end{tikzpicture}
\end{center}
Here $m' = 1$ and $m = 2$, so we must replace $\X^{ext}_1 - \X'^{ext}_1$ with $\pm(\X'^{ext}_2 - \X^{ext}_2)$. We trace the path of dot slides and row operations followed by $\X^{ext}_1$ in {\color{red}{red}} and those followed by $\X'^{ext}_1$ in {\color{blue}{blue}}:

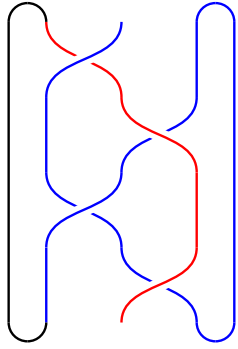
\begin{center}
\begin{tikzpicture}[anchorbase]
	\draw[webs,color=blue] (2,0) to[out=90,in=270] (1,1);
	\draw[line width=5pt,color=white] (1,0) to[out=90,in=270] (2,1);
	\draw[webs,color=red] (1,0) to[out=90,in=270] (2,1);
	\draw[webs,color=blue] (0,0) to (0,1);
	
	\draw[webs,color=blue] (1,1) to[out=90,in=270] (0,2);
	\draw[line width=5pt,color=white] (0,1) to[out=90,in=270] (1,2);
	\draw[webs,color=blue] (0,1) to[out=90,in=270] (1,2);
	\draw[webs,color=red] (2,1) to (2,2);
	
	\draw[webs,color=blue] (1,2) to[out=90,in=270] (2,3);
	\draw[cover] (2,2) to[out=90,in=270] (1,3);
	\draw[webs,color=red] (2,2) to[out=90,in=270] (1,3);
	\draw[webs,color=blue] (0,2) to (0,3);
	
	\draw[webs,color=red] (1,3) to[out=90,in=270] (0,4);
	\draw[cover] (0,3) to[out=90,in=270] (1,4);
	\draw[webs,color=blue] (0,3) to[out=90,in=270] (1,4);
	\draw[webs,color=blue] (2,3) to (2,4);
	
	\draw[webs,color=blue] (2,4) to[out=90,in=180] (2.25,4.25);
	\draw[webs,color=blue] (2.25,4.25) to[out=0,in=90] (2.5,4);;
	\draw[webs,color=blue] (2.5,4) to (2.5,0);
	\draw[webs,color=blue] (2.5,0) to[out=270,in=0] (2.25,-.25);
	\draw[webs,color=blue] (2.25,-.25) to[out=180,in=270] (2,0);
	
	\draw[webs] (0,4) to[out=90,in=0] (-.25,4.25);
	\draw[webs] (-.25,4.25) to[out=180,in=90] (-.5,4);;
	\draw[webs] (-.5,4) to (-.5,0);
	\draw[webs] (-.5,0) to[out=270,in=180] (-.25,-.25);
	\draw[webs] (-.25,-.25) to[out=0,in=270] (0,0);
\end{tikzpicture}
\end{center}

This sequence never applies a row operation using the term $\X^{ext}_1 - \X'^{ext}_1$ which we must replace. This example easily generalizes; given points $\X^{ext}_{m'}$ and $\X'^{ext}_{m'}$, one can trace corresponding paths outward from these points along the component $\L_p$, and these paths will eventually reach the points $\X^{ext}_m$ and $\X'^{ext}_m$.

\end{ex}

\begin{defn} \label{def: red_hhh}
The \textit{reduced colored HOMFLY homology of $\L$ at $p$} is the $\Z_a \times \Z_q \times \Z_t$-graded $R$-module

\[
\overline{HHH}(\L(p)) := H(HH_{\neg p}(F_p(\beta))).
\]
\end{defn}

\begin{rem}
Suppose $\beta$ is non-split and $p$ lies on a strand of minimal color. Then Definition \ref{def: red_hhh} agrees exactly with the reduced HOMFLY homology of \cite{Wed19}. The only challenge is in marrying the two notations: Wedrich's matrix factorization $UMF(\beta)$ is quasi-isomorphic to the Rickard complex $F(\beta)$ (after forgetting the negative differential), Wedrich's $UMF(\beta)/\langle \X'_j = \X_j \rangle$ is exactly our $HH_{\neg p}(F(\beta))$, and Wedrich's $\overline{MF(\L(i))}$ is our $HH_{\neg p}(F_p(\beta))$.
\end{rem}

We pause here to consider the surviving module action on $\overline{HHH}(\L(p))$ before stating our main invariance theorem. With Theorem \ref{thm: point_invariance} in mind, we consider $HH_{\neg p}(F_p(\beta))$ as a complex in $K^b(D^b(\overline{\Scal_\beta(p)}-\Mod{}))$, so $\overline{HHH}(\L(p))$ is a triply-graded $\overline{\Scal_\beta(\L_p)}$-module. The ring $\overline{\Scal_\beta(\L_p)}$ depends heavily on the choice of braid representative $\beta$; however, the induced actions of any two alphabets associated to points on the same component of $\L$ agree, since the actions on $HH_{\neg p}(F_p(\beta))$ complex are homotopic through either dot slides or closures. This leaves an action of the following ring (compare to Definition \ref{def: link_symms}).

\begin{defn}
Given a colored link $\L$ with components colored $(k_1, \dots, k_n)$ and a component $\L_i$ of $\L$, set

\[
\Sym(\X_{\L - \L_i}) := \Sym(\X_1 | \dots | \X_{i - 1} | \X_{i + 1} | \dots | \X_n),
\]
where $\X_j$ is an alphabet of size $k_j$.
\end{defn}

The remainder of Section \ref{subsec: reduced_homology} is dedicated to establishing the following:

\begin{thm} \label{thm: link_invariance}
The reduced HOMFLY homology $\overline{HHH}(\L(p))$ of $\L$ at $p$ is an invariant of the colored link $\L$ and the link component $\L_p \ni p$ up to isomorphism of triply-graded $\Sym(\X_{\L - \L_p})$-modules and overall regrading.
\end{thm}

We adapt the argument of Wedrich to our more general setting. This requires a few technical tools.

\begin{lem} \label{lem: sym_slide}
$f(\X^{ext}) = f(\X'^{ext}) \in W_{\vec{k}}$ for all resolutions $\vec{k}$ of $\beta$ and all symmetric polynomials $f$.
\end{lem}

\begin{proof}
It is easily checked that $f(\X_1 + \X_2) = f(\X'_1 + \X'_2)$ for every ladder web $W_k$ and for every symmetric polynomial $f$ using the relations of Proposition \ref{prop: web_quot_regular} (or Definition \ref{def: new_web_rels}). The result for the whole complex $F(\beta)$ follows from applying this observation at each crossing.
\end{proof}

\begin{cor} \label{cor: Xm_extraneous}
$e_i(\X^{ext}_m) - e_i(\X'^{ext}_m) = 0 \in W_{\vec{k}}/ \langle \X_j - \X'_j \ | \ j \neq m \rangle$ for each $1 \leq i \leq a_m$.
\end{cor}

\begin{proof}
We expand $e_i(\X^{ext})$ as

\[
e_i(\X^{ext}) = e_i(\X^{ext}_m) + \sum_{j = 0}^{i-1} e_j(\X^{ext}_m) e_{i-j} \left( \sum_{k \neq m} \X^{ext}_k \right)
\]
and $e_i(\X'^{ext})$ as

\[
e_i(\X'^{ext}) = e_i(\X'^{ext}_m) + \sum_{j = 0}^{i-1} e_j(\X'^{ext}_m) e_{i-j} \left( \sum_{k \neq m} \X'^{ext}_k \right).
\]
The desired result then follows by inducting on $i$ and canceling the second terms in each expression.

\end{proof}

\begin{prop} \label{prop: reduction_as_cone}
$\prod_{i = 1}^{a_m} (1 + aq^{-2i}) \overline{HHH}(\L(p)) \cong H(K_t(HH(F(\beta)); \X_p))$ as triply-graded $\Sym(\X_{\L-\L_p})$-modules.
\end{prop}

\begin{proof}
Note that $HH(F(\beta)) = K_a(HH_{\neg p} (F(\beta)); \X^{ext}_m - \X'^{ext}_m)$ by definition. By Corollary \ref{cor: Xm_extraneous}, we can apply row operations to the right-hand expression to replace each of the terms $e_i(\X^{ext}_m) - e_i(\X'^{ext}_m)$ in the Koszul differential with $0$, obtaining chain isomorphisms in $K^b(D^b(\Scal_\beta-\Mod{}))$ of the form

\[
HH(F(\beta)) = K_a(HH_{\neg p} (F(\beta)); \X^{ext}_m - \X'^{ext}_m) \cong \prod_{i = 1}^{a_m} (1 + aq^{-2i}) HH_{\neg p}(F(\beta)).
\]
The result then follows from applying the additive functor $H(K_t(-; \X_p))$.
\end{proof}

\begin{lem} \label{lem: fin_gen_module}
$\overline{HHH}(\L(p))$ is a finitely-generated $\Sym(\X_{\L - \L_p})$-module.
\end{lem}

\begin{proof}
It was shown in Lemma 3.45 of \cite{Wed19} that each web $W_{\vec{k}}$ is finitely-generated over $\Scal_\beta$.\footnote{Note that our $\Scal_\beta$ is a quotient of Wedrich's $\Scal_\beta$ by the relations $e_i(\X_1 + \X_2) - e_i(\X'_1 + \X'_2)$ at each crossing. The action on $W_{\vec{k}}$ respects these relations, so this causes us no difficulties.} It follows that $W_{\vec{k}} \tensor[\Scal_{\beta}, der] \overline{\Scal_{\beta}(p)}$ is finitely-generated over $\overline{\Scal_\beta(p)}$. The underlying module of the complex $HH_{\neg p}(F_p(\beta))$ is a finite sum of shifted copies of various $W_{\vec{k}} \tensor[\Scal_\beta, der] \overline{\Scal_\beta(p)}$, so it is also finitely-generated. Since $\Scal_\beta$ is a polynomial ring, it is Noetherian, so passing to homology in the Hochschild and homological directions preserves finite generation.
\end{proof}

Before giving a proof of Theorem \ref{thm: link_invariance}, we cite for later use an extension of Lemma \ref{lem: fin_gen_module}.

\begin{prop} \label{prop: fin_dim}
For any $k$-colored knot $\mathcal{L}$ and marked point $p$, $\overline{HHH}( \mathcal{L}(p))$ is a finite-dimensional $R$-module.
\end{prop}

\begin{proof}
See Proposition 3.46 of \cite{Wed19}.
\end{proof}

\begin{proof}[Proof of Theorem \ref{thm: link_invariance}:]
For a fixed braid word $\beta$ and varying choices of $p$, this is established by Theorem \ref{thm: point_invariance}. Given two pairs $(\beta, p)$ and $(\beta', p')$ such that $\beta$ and $\beta'$ close to the same link, the braids $\beta$ and $\beta'$ must be related by a series of braid relations and Markov moves. Invariance under choice of $p$ for a fixed braid word guarantees that we can take these moves to avoid the point $p$. Then braid relations induce $\Sym(\X_p)$-linear homotopy equivalences of Rouquier complexes, and therefore isomorphisms of reduced homology.

It remains to consider Markov moves avoiding marked points; we adapt the proof from Corollary 3.41 of \cite{Wed19} to our language. Let $(\beta, p)$ and $(\beta', p')$ be two pairs related by such a move. By Theorem \ref{thm: hhh_invariance}, this move induces a $\Sym(\X_\L)$-linear homotopy equivalence $HH(F(\beta)) \simeq HH(F(\beta'))$. This in turn gives an isomorphism

\[
 H(K_t(HH(F(\beta)); \X_p)) \cong H(K_t(HH(F(\beta')); \X_{p'}))
\]
By Proposition \ref{prop: reduction_as_cone}, this induces an isomorphism

\[
\prod_{i = 1}^{a_m} (1 + aq^{-2i}) \overline{HHH}(\hat{\beta}(p))  \cong \prod_{i = 1}^{a_m} (1 + aq^{-2i}) \overline{HHH}(\hat{\beta'}(p')).
\]
Each side is a finitely-generated graded module over the polynomial ring $\Sym(\X_{\L - \L_p})$. The category of such modules is Krull-Schmidt, so we may cancel the factors $\prod_{i = 1}^{a_m} (1 + aq^{-2i})$.
\end{proof}

\subsection{Relation with Unreduced Homology}

A natural question is whether $HHH(\mathcal{L})$ and $\overline{HHH}(\mathcal{L}(p))$ carry the same amount of information. In the uncolored setting, Rasmussen showed in \cite{Ras15} that for all links $\mathcal{L}$ and all points $p$, there is an isomorphism of triply-graded $R$-modules of the form

\begin{equation} \label{eq: ras_isomorphism}
HHH(\mathcal{L}) \cong HHH(\bigcirc) \otimes \overline{HHH}(\mathcal{L}(p)).
\end{equation}

The isomorphism \eqref{eq: ras_isomorphism} has two immediate consequences. First, one can immediately recover the graded dimension of $HHH(\mathcal{L})$ from the graded dimension of $\overline{HHH}(\mathcal{L}(p))$ and vice versa.
Second, since the left-hand side $HHH(\mathcal{L})$ and the factor $HHH(\bigcirc)$ are independent of the choice of component $\mathcal{L}_p \ni p$, so is $\overline{HHH}(\mathcal{L}(p))$.

In this section, we prove an analogous relationship for colored links. Such an isomorphism was originally conjectured by Wedrich in \cite{Wed19}, Remark 3.43. Combining Theorem \ref{thm: main_reduction_thm} with our previous work \cite{Con23, Con24} enables computations of $\overline{HHH}(\mathcal{L}(p))$ for all positive torus knots (and positive torus links with at most one nontrivially colored component); see Section \ref{sec: torus} for details.

\begin{thm} \label{thm: main_reduction_thm}
For any colored link $\L$ and any marked point $p$ on a component $\L_p$ of color $n$, we have an isomorphism of triply-graded $\Sym(\X_{\L - \L_p})$-modules

\[
HHH(\L) \cong HHH(\bigcirc_n) \otimes \overline{HHH}(\L(p)).
\]
\end{thm}

\begin{proof}
We can take the pair $(\beta, p)$ to the following form by a series of Markov moves:

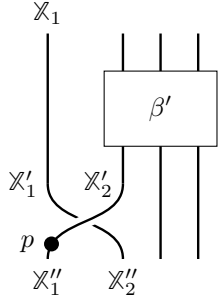
\begin{center}
\begin{tikzpicture}
\draw[webs] (0,0) node[left]{$\X'_2$} to (0,2);
\draw[webs] (.5,-1) to (.5,2);
\draw[webs] (1,-1) to (1,2);
\draw[webs] (-1, 0) node[left]{$\X'_1$} to (-1,2) node[above]{$\X_1$};
\draw[webs] (0,-1) node[below]{$\X''_2$} to[out=90,in=270] (-1,0);
\draw[cover] (-1,-1) to[out=90,in=270] (0,0);
\draw[webs] (-1,-1) node[below]{$\X''_1$} to[out=90,in=270] (0,0);

\node[circle,fill=black,inner sep=2pt,label={west:$p$}] at (-.95,-.8) {};;

\draw [fill=white] (-.25,.5) rectangle++(1.5,1);
\node at (.5,1) {$\beta'$};
\end{tikzpicture}
\end{center}
Here we have included names for various alphabets acting on $F(\beta')$. We denote by $\X$, $\X'$, and $\X''$ the sums of all the alphabets appearing at the same heights as $\X_1$, $\X'_1$, and $\X''_1$, respectively. We will also use the following notation for the alphabets appearing in webs in the Rickard complex associated to the bottom left crossing:

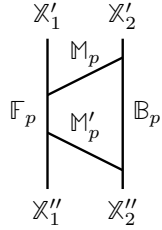
\begin{center}
\begin{tikzpicture}[anchorbase]
\draw[webs] (0,0) node[below]{$\X''_1$}
   -- (0,2) node[pos=0.5,left] {$\F_p$}
   node[above] {$\X'_1$};
\draw[webs] (1,0) node[below]{$\X''_2$} -- (1,2) node[pos=.5, right]{$\B_p$} node[above]{$\X'_2$};
\draw[webs] (1,.25) -- (0,.75) node[pos=.5,above]{$\M'_p$};
\draw[webs] (0,1.25) -- (1,1.75) node[pos=.5,above]{$\M_p$};
\end{tikzpicture}
\end{center}

Since $\X'_1, \X'_2, \X''_1$, and $\X''_2$ all have the same cardinality $n$, we can choose to eliminate the alphabet $\X'_1$ in writing $\Scal_\beta$ (and therefore $\VS_\beta$) as a polynomial algebra, obtaining isomorphisms

\[
\Scal_\beta \cong \Scal_{\beta'} \otimes_R \Sym(\X''_1 | \X''_2); \quad \VS_\beta \cong \VS_{\beta'} \otimes_R \Sym(\X''_1 | \X''_2 | \B_p).
\]
Then $\overline{\Scal_\beta(p)} \cong \Scal_{\beta'} \otimes_R \Sym(\X''_2)$. As usual, $\Scal_\beta \cong \overline{\Scal_\beta(p)} \otimes_R \Sym(\X''_1)$, so we have inclusion and projection maps $\iota \colon \overline{\Scal_\beta(p)} \hookrightarrow \Scal_\beta$ and $\pi \colon \Scal_\beta \to \overline{\Scal_\beta(p)}$. We denote by $i_!$ and $\pi_!$ the induced extension of scalar functors on the (homotopy category of the) derived category:

\[
i_! := - \otimes_{\overline{\Scal_\beta(p)}} \Scal_\beta = - \otimes_R \Sym(\X''_1); \quad \pi_! := - \tensor[\Scal_\beta, der] \overline{\Scal_\beta(p)}.
\]

In this notation, the reduced HOMFLY homology of $\L$ at $p$ is exactly

\[
\overline{HHH}(\L(p)) := H(HH_{\neg p}(\pi_!(F(\beta)))).
\]
Set $\overline{\VS_\beta(p)} := \pi_!(\overline{\Scal_\beta(p)})$; then $\overline{\VS_\beta(p)} \cong \VS_{\beta'} \otimes_R \Sym(\X''_2 | \B_p)$, and $\VS_\beta \cong \overline{\VS_\beta(p)} \otimes_R \Sym(\X''_1)$. For each web $W_{\vec{k}}$ appearing in the Rickard resolution of $\beta$, we denote by $k_p$ the highest degree of the surviving generators \footnote{Or just the strand label in the traditional web description of Proposition \ref{prop: web_quot_regular}.} $e_i(\B_p)$ and by $W_{\vec{k} - k_p}$ the corresponding web appearing in the Rickard resolution of $\beta'$. We point out that $W_{\vec{k} - k_p} \otimes_R \Sym(\X''_1 | \X''_2 | \B_p)$ is isomorphic to the quotient of $\VS_\beta$ by (the ideal generated by) all the generators of $\JS_{\vec{k}}$ which did not arise from the ideal $\JS_{k_p}$ associated to the bottom-left crossing of $\beta$.

For each choice of resolution $\vec{k}$, almost all of the generators of $\JS_{\vec{k}}$ lie in $\overline{\VS_\beta(p)}$. The only outlier is the family of generators

\begin{equation} \label{eq: problematic_gens}
\{e_i(\X'_1 + \X'_2) - e_i(\X''_1 + \X''_2) \ | \ n < i \leq k_p \}.
\end{equation}

Within $\VS_\beta$, we already have $e_i(\X'_1 + \X'_2) = e_i(\X''_1 + \X''_2)$ for all $1 \leq i \leq n$. Since $\X'_j = \X''_j$ for all $j \neq 1, 2$, we also have $e_i(\X') = e_i(\X'') \in \VS_\beta$. After passing to the quotient $W_{\vec{k} - k_p} \otimes_R \Sym(\X''_1 | \X''_2 | \B_p)$, we have $e_i(\X) = e_i(\X') = e_i(\X'')$ by Lemma \ref{lem: sym_slide} applied to $\beta'$. In the further quotient by $\langle \X_j - \X''_j \ | \ j \neq 1 \rangle$, $e_j(\X - \X_1) = e_j(\X'' - \X''_1)$ for each $j$. In total, we obtain

\begin{align*}
e_i(\X) & = \sum_{j + k = i} e_j(\X_1) e_k(\X - \X_1) \\
& = \sum_{j + k = i} e_j(\X'_1) e_k(\X'' - \X''_1) \\
= e_i(\X'') & = \sum_{j + k = i} e_j(\X''_1) e_k(\X'' - \X''_1) \in \left( W_{\vec{k} - k_p} \otimes_R \Sym(\X''_1 | \X''_2 | \B_p) \right) /\langle \X_j - \X''_j \ | \ j \neq 1 \rangle.
\end{align*}

Canceling the terms $e_k(\X'' - \X''_1)$ on the left and right, we see that $e_j(\X_1') = e_j(\X''_1)$ for all $1 \leq j \leq n$ in this quotient. By the relation $e_i(\X'_1 + \X'_2) = e_j(\X''_1 + \X''_2)$ for $1 \leq i \leq n$, this forces $e_j(\X'_2) = e_j(\X''_2)$ as well. As a consequence, each generator from \eqref{eq: problematic_gens} vanishes upon passing to the quotient $\left(W_{\vec{k} - k_p} \otimes_R \Sym(\X''_1 | \X''_2 | \B_p) \right)/\langle \X_j - \X''_j \ | \ j \neq 1 \rangle$.

Now, each generator of the denominator of this quotient corresponds to a term not arising from $\JS_{k_p}$ in the Koszul differential of the complex

\[
HH_{\neg p}(W_{\vec{k}}) \cong K_a(\VS_{\beta}; \JS_{\vec{k}} \sqcup \{\X_j - \X''_j \ | \ j \neq 1\}).
\]
Consequently, we can apply row operations within this complex to replace the terms \eqref{eq: problematic_gens} in the differential with $0$. All of the remaining terms in the differential lie in $\overline{\VS_\beta(p)}$, so they are preserved by applying the composition $i_! \circ \pi_!$. Since $\VS_\beta \cong \overline{\VS_\beta(p)} \otimes_R \Sym(\X''_1)$ is also preserved by $i_! \circ \pi_!$, we obtain

\[
i_! \circ \pi_! \left( HH_{\neg p}(W_{\vec{k}}) \right) \cong HH_{\neg p}(W_{\vec{k}}).
\]
Applying these quasi-isomorphisms in each homological degree, we obtain a homotopy equivalence

\[
i_! \circ \pi_! \left( HH_{\neg p}(F(\beta)) \right) \simeq HH_{\neg p}(F(\beta)).
\]
By the proof of Proposition \ref{prop: reduction_as_cone}, this induces an isomorphism

\[
HHH(\L) = H(HH(F(\beta))) \cong \prod_{i = 1}^n (1 + aq^{-2i}) H\left(  i_! \circ \pi_! \left( HH_{\neg p}(F(\beta)) \right) \right)
\]

Now, since $\Scal_\beta$ is a flat $\overline{\Scal_\beta(p)}$-module, $i_!$ is exact, and therefore commutes with taking homology. Further, $HH_{\neg p}$ is just passing to a Koszul complex, so it commutes with $\pi_!$. We can therefore rewrite the above isomorphism as

\begin{align*}
HHH(\L) & \cong \prod_{i = 1}^n (1 + aq^{-2i}) i_! \left( H \left( HH_{\neg p} \left( \pi_! (F(\beta)) \right) \right) \right) \\
& \cong \prod_{i = 1}^n (1 + aq^{-2i}) \overline{HHH}(\L(p)) \otimes_R \Sym(\X''_1) \\
& \cong \prod_{i = 1}^n \left( \cfrac{1 + aq^{-2i}}{1 - q^{2i}} \right) \overline{HHH}(\L(p)).
\end{align*}

Finally, since all of the above isomorphisms were $\Sym(\X_{\L - \L_p})$-equivariant, this is an isomorphism of triply-graded $\Sym(\X_{\L - \L_p})$-modules.
\end{proof}

As in the uncolored case, we have the following immediate corollary:

\begin{cor}
$\overline{HHH}(\mathcal{L}(p)) \cong \overline{HHH}(\mathcal{L}(p'))$ as triply-graded $R$-modules for any two marked points $p, p'$ lying on components of $\mathcal{L}$ of the same color.
\end{cor}

\section{Torus Links} \label{sec: torus}

\subsection{Reduced Colored Torus Knot Homology}

Colored HOMFLY homology comes in many different flavors. The invariants described in Section \ref{sec: link_hom} are \textit{intrinsically colored} homologies, defined using Rickard complexes of singular Soergel bimodules. There are also \textit{projector-colored} homologies, defined by replacing a $k$-colored strand with a $k$-fold cabling of $1$-colored strands, then inserting a certain complex called a \textit{categorified projector}; see \cite{Hog18, AH17, Cau17, EH17b, Con24}. Both of these types admit a deformation called a \textit{$y$-ification}; see \cite{GH22, HRW24} for this deformation in the intrinsically colored setting and \cite{Con23, BGHW25} for the projector-colored setting. In \cite{GW19}, Gorsky--Wedrich construct what one might call \textit{invariant colored} homology, in which the coloring is accounted for by taking invariants under an action of the symmetric group on HOMFLY homology of cables. Exploring the landscape of different invariants and their relations presents an interesting problem. At the level of finite rank $\mathfrak{gl}(N)$-homologies, the corresponding problem is largely solved thanks to the toolkit of categorical representation theory and skew Howe duality; see e.g. \cite{MW18, BHPW23, Mer25} for discussions.

We established in \cite{Con24} isomorphisms among various versions of projector-colored and intrinsically colored HOMFLY homologies, including finite and infinite versions of the former and deformations of both types, allowing one to determine each of these invariants for a given link from any of the others. Theorem \ref{thm: main_reduction_thm} allows us to add reduced HOMFLY homology to this list.

In particular, we directly computed in \cite{Con23} the $y$-ified, infinite projector-colored HOMFLY homology of all positive torus knots. Armed with Theorem \ref{thm: main_reduction_thm}, we can use this chain of isomorphisms to compute the reduced homology of these knots as well.

\begin{notn} \label{notn: gr_dim}
For each $m, n \geq 0$, let $d := \gcd(m, n)$. We denote by $T(m, n)_{k_1, \dots, k_d}$ the positive torus link $T(m, n)$ with components colored $k_1, \dots, k_d$, and we set

\[
\mathcal{P}(m, n)_{k_1, \dots, k_d}(A, Q, T) := \dim \Big( \overline{HHH} \big( T(m, n)_{k_1, \dots, k_d}(p) \big) \Big)
\]
where we always assume that $p$ lies on the $k_1$-colored component of $T(m, n)$.
\end{notn}

\begin{thm} \label{thm: red_torus_links}
For each pair $v \in \{0, 1\}^{l + m}$ and $w \in \{0, 1\}^{l + n}$ with $|v| = |w| = l$ and each $\sigma \in \SG_l$, let $p_{\sigma}(v, w) \in \mathbb{N}[[A,Q,T]]$ be the unique power series satisfying the following recursions:

\begin{enumerate}
\item $p_e(0^m, \emptyset) = \left( \cfrac{1 + A}{(1 - Q)(1 - T)} \right)^m$ and $p_e(\emptyset, 0^n) = \left( \cfrac{1 + A}{(1 - Q)(1 - T)} \right)^n$;

\item $p_e(0^m1, 0^n1) = \left( \cfrac{1 + A}{(1-Q)(1-T)} \right) p_e(0^m, 0^n)$ for all $m, n \geq 0$;

\item $p_\sigma(v1,w1) = \left( \cfrac{Q^l + A}{1 - Q} \right) p_{Tr(\sigma)}(v, w)$ if $l \geq 1$ and $\sigma(n) = n$;

\item $p_\sigma(v1,w1) = (Q^l + A) p_{Tr(\sigma)}(v, w)$ if $l \geq 1$ and $\sigma(n) \neq n$;

\item $p_\sigma(v0, w1) = p_{\sigma \pi_{(10^{l - 1})}}(v, 1w)$ if $l \geq 1$;

\item $p_\sigma(v1,w0) = p_{\pi_{(10^{l-1})}\sigma} (1v, w)$ if $l \geq 1$;

\item $p_e(0^m, 0^n) = p_e(10^{m - 1}, 10^{n - 1})$ for all $m, n \geq 1$;

\item $p_\sigma(v0, w0) = Q^{-l} p_{e_1 \sqcup \sigma}(1v, 1w) + Q^{-l} T p_\sigma(0v, 0w)$ if $l \geq 1$.
\end{enumerate}

For each $k \geq 0$, set

\[
f(k) := [k]! \left( \prod_{i = 1}^k \cfrac{1 + AQ^{1-i}}{1-Q^i}\right).
\]

Further, for each $k, m, n \geq 1$, set $d = \gcd(m, n)$ and

\[
v_k := 1^k0^{md^{-1}(d - 1) + k(md^{-1} - 1)}; \quad w_k := 1^k0^{nd^{-1}(d - 1) + k(nd^{-1} - 1)}.
\]

Then for all $k, m, n$, we have

\[
\mathcal{P}(m, n)_{k, 1, \dots, 1} = (1-T)^d f(k)^{-1} p_e(v_k, w_k).
\]
\end{thm}

\begin{rem}
For experts, we briefly explain the multiplicative correction terms in Theorem \ref{thm: red_torus_links}. The factor of $[k]!$ from $f(k)$ comes from passing from projector-colored to intrinsically colored homology, and $\left( \prod_{i = 1}^k \cfrac{1 + AQ^{1-i}}{1-Q^{2i}}\right)$ comes from dividing by the dimension of the colored homology of the unknot as in Theorem \ref{thm: main_reduction_thm}. The factors of $1 - T$ have two sources. For uncolored components, this exactly corrects for the dimension of the polynomial ring $R[y]$ in the $y$-ified complex; this effect contributes a factor of $(1-T)^{d-1}$ in total. The remaining factor of $1 - T$ comes from replacing a $1$-colored strand with a \textit{renormalized} finite column projector on a single strand (c.f. Remark 4.14 in \cite{Con24}). 
\end{rem}

\subsection{Doubly-Graded Exponential Growth}

This section is dedicated to the proof of our second main theorem:

\begin{thm} \label{thm: doubly_graded_growth}
For any $m, n \geq 1$ with $\gcd(m, n) = 1$, the reduced HOMFLY homology of the $k$-colored positive torus knot satisfies

\begin{align} \label{eq: doubly_graded_growth}
\mathcal{P}(m, n)_k \Big|_{Q = 1} = \left( \mathcal{P}(m, n)_1 \Big|_{Q = 1} \right)^k.
\end{align}
\end{thm}

We will prove Theorem \ref{thm: doubly_graded_growth} by setting $Q = 1$ directly in the recursion of Theorem \ref{thm: red_torus_links}. The terms $1 - Q$ in the denominators of rules (1), (2), and (3) present some difficulties which we have to overcome. This requires us to discuss the permutations appearing in this recursion. We will keep the discussion brief; see \cite{Con23} for more details.

The functions $p_{\sigma}(v, w)$ arise as the dimensions of a $y$-ification of the Hochschild homology of certain complexes $\textbf{C}^y(v, w)_{\sigma}$ from \cite{Con23}. The complexes are built as follows:

\begin{center}
\begin{gather*}
\textbf{C}^y(v, w)_{\sigma} := \quad
\begin{tikzpicture}[baseline=(current bounding box.center),scale=.75]
	\draw[webs] (1,0) node[below]{$n$} to (1,1);
	\draw[webs] (1,1) to[out=90,in=270] (0,3);
	\draw[webs] (0,0) node[below]{$m$} to (0,1);
	\draw[color=white, line width=5pt] (0,1) to[out=90,in=270] (1,3);
	\draw[webs] (0,1) to[out=90,in=270] (1,3);
	\draw[webs] (2,0) node[below]{$l$} to (2,1);
	\node[draw, fill=white, minimum width=1cm] at (1.5,.5) {$\beta_w$};
	\draw[webs] (2,1) to node[pos=.5, draw, fill=white]{$\hat{K}^y_{l, \sigma}$} (2,3);
	\draw[webs] (0,3) to (0,4) node[above]{$n$};
	\draw[webs] (1,3) to (1,4) node[above]{$m$};
	\draw[webs] (2,3) to (2,4) node[above]{$l$};
	\node[draw, fill=white, minimum width=1cm] at (1.5,3.5) {$\alpha_v$};
\end{tikzpicture}
\end{gather*}
\end{center}

In this context, the labels $m, n, l$ represent (blackboard framed) cables of $1$-colored strands and \textit{not} colored strands\footnote{One can remember this by observing that the resulting complex would not be balanced, since $(n, m, l) \neq (m, n, l)$ in general.}. $\alpha_v$ and $\beta_w$ are \textit{shuffle braids} on $m + l$ and $n + l$ strands, respectively, which are positive braid lifts of certain shuffle permutations $\pi_v$ and $\pi_w$ which can be read off from $v$ and $w$. Under this dictionary, the digits in $v$ and $w$ correspond exactly to strands entering $\alpha_v$ and $\beta_w$ from the top and bottom, respectively, read from left to right. The complex $\hat{K}^y_{l, \sigma}$ is the \textit{renormalized twisted finite $y$-ified column projector} of \cite{Con23}. Its exact form is not important for our purposes; what is important is that it has an associated permutation $\sigma$.

Part of the definition of a $y$-ified complex on $n$ strands is an associated permutation $\omega \in \SG_n$, which we can often view as an explicit permutation taking the strands on the bottom to the strands on the top. We can assign to each triple of input data the associated permutation $\omega(v, w, \sigma) \in \SG_{m + l + n}$ underlying $\textbf{C}^y(v, w)_{\sigma}$. This allows us to group the strands of $\textbf{C}^y(v, w)_{\sigma}$, and therefore the digits of $v$ and of $w$, by the orbits of $\omega(v, w, \sigma)$ in which they lie.

Each of the rules (1)-(8) from Theorem \ref{thm: red_torus_links} corresponds to some explicit manipulation of the complex $\textbf{C}^y(v, w)_{\sigma}$, and we can track how $\omega(v, w, \sigma)$ evolves as we perform these manipulations. The following properties of this evolution will be crucial; they can be read off directly from the description of \cite{Con23}.

\begin{lem} \label{lem: recursion_facts} The permutations $\omega(v, w, \sigma)$ satisfy:

\begin{enumerate}
\item $\omega(v_k, w_k, e)$ has exactly $dk$ orbits for all $k, m, n$, where here $d = \gcd(m, n)$.

\item Rules (1), (2), and (3) from the recursion of Theorem \ref{thm: red_torus_links} are applied \textit{exactly} when the number of orbits of $\omega(v, w, \sigma)$ decreases by $1$.

\item The final digits of $v$ and $w$ always correspond to strands which lie in the same orbit of $\omega(v, w, \sigma)$.

\item All digits which are not involved in a given recursion rule remain in the same orbit before after applying the rule.

\item When a digit of $v$ or $w$ is moved from the end of its string to the beginning (using rules (5), (6), (7), and (8)), the orbit containing the strand which corresponds to that digit does not change.
\end{enumerate}
\end{lem}

We are now prepared to state a simplified form of our recursion which directly computes $\mathcal{P}(m, n)_k \Big|_{Q = 1}$.

\begin{prop} \label{prop: simpler_recursion}

For each $v, w$ as in Theorem \ref{thm: red_torus_links}, let $\p(v, w) \in \mathbb{N}[A, T]$ be the unique polynomial satisfying the following recursions:

\begin{enumerate}
\item $\p(0^m, \emptyset) = (1+A)^m$ and $\p(\emptyset, 0^n) = (1+A)^n$;

\item $\p(v1, w1) = (1 + A) \p(v, w)$;

\item $\p(v0, w1) = \p(v, 1w)$;

\item $\p(v1, w0) = \p(1v, w)$;

\item $\p(0^m, 0^n) = \p(10^{m - 1}, 10^{n - 1})$ for $m, n \geq 1$;

\item $\p(v0, w0) = \p(1v, 1w) + T\p(0v, 0w)$ if $l \geq 1$.
\end{enumerate}

Then for each $k, m, n \geq 1$ with $\gcd(m, n) = 1$, we have

\[
\mathcal{P}(m, n)_k \Big|_{Q = 1} = (1 + A)^{-k} \p(v_k, w_k).
\]

\end{prop}

\begin{proof}
We begin by dealing with the factor $f(k)^{-1}$ from Theorem \ref{thm: red_torus_links}. We can simplify $f(k)$ by factoring the denominator, obtaining

\begin{align*}
f(k) & = [k]! \left( \prod_{i = 1}^k \cfrac{1 + AQ^{1 - i}}{1 - Q^i} \right) \\
& = \left( \prod_{i = 1}^k \cfrac{[i]}{1 + Q + \dots + Q^{i - 1}} \right) \left( \prod_{i = 1}^k \cfrac{1 + AQ^{1 - i}}{1 - Q} \right) \\
& = Q^{\bullet} (1-Q)^k \left( \prod_{i = 1}^k 1 + AQ^{1 - i} \right)
\end{align*}
up to an overall shift of the form $Q^\bullet$. The first two properties of Lemma \ref{lem: recursion_facts} guarantee that the power of $1 - Q$ appearing in the denominator of each summand of $p_e(v_k, w_k)$ is exactly $k$, so we are free to cancel these factors against the factor of $(1 - Q)^k$ in $f(k)$.

After this cancellation, we can set $Q = 1$ directly in Theorem \ref{thm: red_torus_links}. This transforms the correction term $(1 - T)f(k)^{-1}$ of Theorem \ref{thm: red_torus_links} to $(1 - T)(1 + A)^{-k}$. This also collapses Rules (3) and (4) of Theorem \ref{thm: red_torus_links}, so that both now read $p_{\sigma}(v1, w1) = (1+A) p_{Tr(\sigma)}(v, w)$. Since these two rules were the only instance of the recursion depending on the permutation $\sigma$, we're now free to drop $\sigma$ entirely. Rules (5), (6), (7), and (8) of Theorem \ref{thm: red_torus_links} now become the new rules (4), (5), (6), and (7), respectively.

To conclude, we'd like to cancel the factor of $1 - T$ in the new correction terms against the remaining factor of $1 - T$ from the denominators of rules (1) and (2) in Theorem \ref{thm: red_torus_links}. Indeed, after multiplying by $f(k)^{-1}$ and setting $Q = 1$, we can group the summands of $p_e(v_k, w_k)$ according the power of $1 - T$ appearing in their denominators:

\[
\left(f(k)^{-1} p_e(v_k, w_k) \right) \Big|_{Q = 1} = \sum_{r \geq 0} \cfrac{p_e(v_k, w_k)_r}{(1 - T)^r}.
\]
The recursion of Theorem \ref{thm: red_torus_links} can only end by applying rule (1) or (2), so we must have $p_e(v_k, w_k)_0 = 0$. Meanwhile, after multiplying by the correction term $1 - T$, we obtain

\[
\mathcal{P}(m, n)_k \Big|_{Q = 1} = \sum_{r \geq 1} \cfrac{p_e(v, w)_r}{(1 - T)^{r - 1}}.
\]

By Proposition \ref{prop: fin_dim}, $\overline{HHH}(T(m, n)_k(p))$ is always finite-dimensional, so $\mathcal{P}(m, n)_k \Big|_{Q = 1}$ is a polynomial in $A$ and $T$ rather than a formal power series. The summands with index $r \geq 2$ in the above sum must therefore cancel. Since each numerator $p_e(v_k, w_k)_r$ has only positive coefficients, the only way this cancellation can happen is if $p_e(v_k, w_k)_r = 0$ for all $r \geq 2$. Then every summand of $\left( f(k)^{-1} p_e(v, w) \right) |_{Q = 1}$  must have a denominator of exactly $1 - T$, so we are free to ignore this factor in both the recursion and the correction terms.
\end{proof}

Recall that for each pair of sequences $v, w$ and choice of permutation $\sigma \in \SG_l$, we can partition the digits of $v$ and of $w$ according to the orbits of $\omega(v, w, \sigma)$. We will denote by $v^i_\sigma$ and $w^i_\sigma$ the subsequences of $v$ and of $w$ lying in the $i^{th}$ such orbit. This partition depends on the choice of permutation $\sigma$.

\begin{lem} \label{lem: rec_mult}
For each $\sigma \in \SG_l$, if the recursions which compute $\p(v, w)$ and each $\p(v^i_\sigma, w^i_\sigma)$ never use Rule (5) of Proposition \ref{prop: simpler_recursion}, then $\p(v, w) = \prod_i \p(v^i_\sigma, w^i_\sigma)$. 
\end{lem}

\begin{proof}
We consider the effect of each of the rules from Proposition \ref{prop: simpler_recursion} on the product $\prod_i \p(v^i_\sigma, w^i_\sigma)$, beginning with Rule (1).

\textbf{Rule (1):} If $w$ is empty (so $v = 0^m$), then so are the subsequences $w^i$ for each $i$, and each subsequence of $v$ is of the form $v^i = 0^{n_i}$. It follows that $\p(v, \emptyset) = (1 + A)^m$, while $\p(v^i, \emptyset) = (1 + A)^{n_i}$. Now observe that $\sum_i n_i = m$, so the desired multiplicativity holds in this case.

\vspace{1em}

\textbf{Rules (2), (3), (4), and (6):} Each of these rules replaces a pair $(v, w)$ with a linear combination of new pairs $(v', w')$, and this modification depends only on the final digits of $v$ and $w$. To define new subsequences $(v')^i_{\sigma'}, (w')^i_{\sigma'}$, we need to choose a new permutation $\sigma'$. We do so according to the prescription in Theorem \ref{thm: red_torus_links}; this ensures that the properties of Lemma \ref{lem: recursion_facts} apply. With this standard in place, we suppress the permutations, writing just $v^i$ and $w^i$, for the remainder of the proof.

In each step, the final digits of $v$ and $w$ lie in the same orbits $v^i$, $w^i$ by Property (3) of Lemma \ref{lem: recursion_facts}. If these digits are recycled to the beginning of the sequences to form $v'$ and $w'$ (as in Rules (3), (4), and (6)), the recycled digits must then lie in $(v')^i$ and $(w')^i$ by Property (5) of Lemma \ref{lem: recursion_facts}. By Property (4) of Lemma \ref{lem: recursion_facts}, the other subsequences are unchanged; that is, $(v^j, w^j) = ((v')^j, (w')^j)$ for all $j \neq i$.

To recap, each of these rules transforms only a single pair of subsequences $(v^i, w^i)$, and does so in a way which exactly reproduces the transformation of the pair $(v, w)$. We can therefore view these rules as applying independently to the relevant factors $p(v^i, w^i)$ in the product and leaving the other factors fixed. 
\end{proof}

\begin{proof}[Proof of Theorem \ref{thm: doubly_graded_growth}]:
By Proposition \ref{prop: simpler_recursion}, it suffices to show 

\[
\p(v_k, w_k) = \p(v_1, w_1)^k.
\]
An easy analysis of the complexes $\textbf{C}^y(v_k, w_k)_{e}$ shows that $(v_k)^i_e = (10^{m - 1}) = v_1$ and $(w_k)^i_e = (10^{n - 1}) = w_1$ for each $i$. If the recursion which computes $\p(v_k, w_k)$ never uses Rule (5), we can therefore conclude by applying Lemma \ref{lem: rec_mult}.

It remains to show that the recursion which computes $\p(v_k, w_k)$ never uses Rule (5). Equivalently, we wish to show that no intermediate pair of nonempty sequences $(\tilde{v}, \tilde{w})$ appearing in this recursion consist entirely of $0$'s. As shown in the proof of Lemma \ref{lem: rec_mult}, the other rules all apply locally to the subsequences $v^i$ and $w^i$, so it suffices to show that no nonempty pair of subsequences $(\tilde{v}^i, \tilde{w}^i)$ consist entirely of $0$'s. Since $(v_k^i, w_k^i) = (v_1, w_1)$, this reduces to checking that we never use Rule (5) in evaluating $\p(v_1, w_1)$.

If we did use Rule (5) in evaluating $\p(v_1, w_1)$, then in the step before it was applied, we would have used Rule (2) to evaluate a term of the form $\p(0^r1, 0^s1)$. The only choice of permutation $\sigma$ here is the trivial one. Making this choice, the application of Rule (2) would then decrease the number of orbits. Since $\omega(v_1, w_1, e)$ started with only one orbit and no rule increases the number of orbits, we must be left with an empty sequence. Hence Rule (5) is never used.
\end{proof}

\subsubsection{$T(2, 2n+1)$}

For torus knots $T(2, 2n+1)$, we can implement the recursion of Proposition \ref{prop: simpler_recursion} to obtain a closed form expression for the doubly-graded dimension of $\overline{HHH}$ after setting $Q = 1$.

\begin{thm}
For each $n, k \geq 1$, we have

\[
\mathcal{P}(2, 2n + 1)_k \Big|_{Q = 1} = \left(T^n + (1 + A) \sum_{j = 0}^{n - 1}T^j \right)^k.
\]

\end{thm}

\begin{proof}
By Theorem \ref{thm: doubly_graded_growth}, it suffices to show that $\mathcal{P}(2, 2n + 1)_1 \Big|_{Q = 1} = T^n + (1 + A) \sum_{j = 0}^{n-1} T^j$. By Proposition \ref{prop: simpler_recursion}, we have

\[
\mathcal{P}(2, 2n + 1)_1 \Big|_{Q = 1} = (1 + A)^{-1} \p(10, 10^{2n}).
\]

We use the recursion of Proposition \ref{prop: simpler_recursion} to explicitly evaluate $\p(10, 10^{2n})$. We can simplify any expression of the form $\p(10, w00)$ by applying first Rule (6), then Rule (4), obtaining

\begin{equation} \label{eq: rec_simp}
\p(10, w00) = \p(11, 1w0) + T \p(01, 0w0) = \p(11, 1w) + T \p(10, 0w).
\end{equation}
In particular, we can apply \eqref{eq: rec_simp} when $w = 10^{2(n - 1)}$, obtaining

\[
\p(10, 10^{2n}) = \p(11, 110^{2(n - 1)}) + T \p(10, 010^{2(n - 1)}).
\]

Now if $n > 1$, we can once again apply \eqref{eq: rec_simp} to the second summand by taking $w = 10^{2(n - 2)}$. This process continues until all the trailing $0$'s of $10^{2n}$ have been exhausted, leaving

\[
\p(10, 10^{2n}) = \left(\sum_{j = 0}^{n - 1} T^j \p(11, 10^j10^{2(n - 1 - j)}) \right) + T^n \p(10, 0^n1).
\]
By repeatedly applying Rules (2) and (4), we eventually see $\p(11, w) = (1 + A)^2$ for any $w$. Meanwhile, by Rule (3), we have $\p(10, 0^n1) = \p(1, 10^n)$, which evaluates to $1 + A$ after $n$ applications of Rule (4) and $1$ application of Rule (2). In total, we see

\[
\p(10, 10^{2n}) = \left(\sum_{j = 0}^{n - 1} T^j (1 + A)^2 \right) + T^n (1 + A) = (1+A) \left( T^n + (1 + A) \sum_{j = 0}^{n - 1} T^j \right).
\]
Multiplying by $(1 + A)^{-1}$ then gives exactly the desired result.

\end{proof}

\begin{rem}
In \S 7.2 of \cite{GGS18}, GGS give a conjectural closed form for $\mathcal{P}(2, 2n+1)_k \Big|_{A = 0}$ as a polynomial in $Q$ and $T$ involving $Q$-binomial coefficients. It would be interesting to try to extract their expression directly from the recursion of Theorem \ref{thm: red_torus_links}.
\end{rem}

\subsection{Color Shifting}

Given an $d$-component torus link $T(dm, dn)$ with $\gcd(m, n) = 1$, recall that each component of $T(dm, dn)$ is isotopic to $T(m, n)$. In particular, if $m = 1$ or $n = 1$, then each component is isotopic to the unknot, and Conjecture \ref{conj: color_shift} should hold. When only one component of such a link is nontrivially colored, we can test this conjecture against the results of Theorem \ref{thm: red_torus_links}. Recall that we denote by $\mathcal{P}(dm, dn)_{k_1, \dots, k_d}$ the graded dimension of the reduced HOMFLY homology of the positive torus link $T(dm, dn)$ with link components colored $k_1, \dots, k_d$, where the chosen point $p$ always lies on the $k_1$-colored component.

\subsubsection{$T(2, 2n)$}

We begin by testing $2$-component torus links of the form $T(2, 2n)$. When $n = 1$, direct computation gives the following:

\begin{align*}
\mathcal{P}(2,2)_{1,1} & = \cfrac{Q(1 - T) + A + T}{1 - Q}; \\
\mathcal{P}(2,2)_{2,1} & = \cfrac{Q^2(1 - T) + A + T}{1 - Q}; \\
\mathcal{P}(2,2)_{3,1} & = \cfrac{Q^3(1 - T) + A + T}{1 - Q}; \\
\mathcal{P}(2,2)_{4,1} & = \cfrac{Q^4(1 - T) + A + T}{1 - Q}.
\end{align*}
This pattern suggests that $T(2,2)$ satisfies

\begin{equation} \label{eq: T22_shifting}
\mathcal{P}(2,2)_{k,1} = \cfrac{Q^k(1 - T) + A + T}{1 - Q}.
\end{equation}
Notice that when $k = 1$, we obtain $\mathcal{P}(2,2)_{0,1} = \cfrac{1 + A}{1 - Q}$, which is exactly the unreduced HOMFLY homology of the $1$-colored unknot. We have verified \eqref{eq: T22_shifting} for all $k \leq 100$. 

When $n = 2$, we obtain

\begin{align*}
\mathcal{P}(2,4)_{1,1} & = \cfrac{Q^2(1 - T) + Q(A + T)(1 - T) + T(A + T)}{1 - Q}; \\
\mathcal{P}(2,4)_{2,1} & = \cfrac{Q^4(1 - T) + Q^2(A + T)(1 - T) + T(A + T)}{1 - Q}; \\
\mathcal{P}(2,4)_{3,1} & = \cfrac{Q^6(1 - T) + Q^3(A + T)(1 - T) + T(A + T)}{1 - Q}; \\
\end{align*}
This pattern suggests that $T(2,4)$ satisfies

\begin{equation} \label{eq: T24_shifting}
\mathcal{P}(2,4)_{k,1} = \cfrac{Q^{2k}(1 - T) + Q^k(A + T)(1 - T) + T(A + T)}{1 - Q}.
\end{equation}
Again, when $k = 0$, we recover exactly $\mathcal{P}(2,4)_{0,1} = \cfrac{1+A}{1-Q}$. We have verified \eqref{eq: T24_shifting} for all $k \leq 90$.

Once more, when $n = 3$, we obtain

\begin{align*}
\mathcal{P}(2,6)_{1,1} = \cfrac{Q^3(1 - T) + Q^2 (1 - T)(A + T) + QT (1 - T)(A + T) + T^2 (A + T)}{1-Q}; \\
\mathcal{P}(2,6)_{2,1} = \cfrac{Q^6(1 - T) + Q^4 (1 - T)(A + T) + Q^2T (1 - T)(A + T) + T^2 (A + T)}{1-Q}; \\
\mathcal{P}(2,6)_{3,1} = \cfrac{Q^9(1 - T) + Q^6 (1 - T)(A + T) + Q^3 T (1 - T)(A + T) + T^2 (A + T)}{1-Q}. \\
\end{align*}
In general, this suggests that $T(2, 6)$ satisfies

\begin{equation} \label{eq: T26_shifting}
\mathcal{P}(2, 6)_{k,1} = \cfrac{Q^{3k}(1 - T) + Q^{2k} (1 - T)(A + T) + Q^k T (1 - T)(A + T) + T^2 (A + T)}{1-Q}.
\end{equation}
We have verified \eqref{eq: T26_shifting} for all $k \leq 90$.

The above results also exhibit some regularity in $n$, suggesting that $k$-colored $T(2, 2n)$ satisfies

\begin{equation} \label{eq: T2N_shifting}
\mathcal{P}(2, 2n)_{k, 1} = \cfrac{Q^{nk} (1 - T) + \sum_{j = 1}^{n - 1} Q^{jk} T^{n - j - 1} (1 - T)(A + T) + T^{n - 1} (A + T)}{1 - Q}.
\end{equation}
Plugging $k = 0$ into \eqref{eq: T2N_shifting} indeed recovers $\mathcal{P}(2, 2n)_{0, 1} = \cfrac{1+A}{1-Q}$. We have verified \eqref{eq: T2N_shifting} for all $n \leq 15$ and $k \leq 20$.

\subsubsection{$T(3, 3)$}

Again by direct computation, we obtain

\begin{align*}
\mathcal{P}(3,3)_{1,1,1} & = \cfrac{Q^3(1 - T)^2 + Q^2(A + T)(1 - T)^2 + Q((A + T)(1 - T^2) + T(A + T)(1 - T)) + (A + T^2)(A + T)}{(1 - Q)^2}; \\
\mathcal{P}(3,3)_{2,1,1} & = \cfrac{Q^5(1 - T)^2 + Q^3(A + T)(1 - T)^2 + Q^2(A + T)(1 - T^2) + QT(A + T)(1 - T) + (A + T^2)(A + T)}{(1 - Q)^2}. \\
\end{align*}
This suggests that $T(3, 3)$ satisfies

\begin{multline} \label{eq: T33_shifting}
(1-Q)^2 \mathcal{P}(3, 3)_{k,1,1} = Q^{2k}(Q(1 - T)^2) + Q^k(Q(T+A)(1-T)^2 + (A+T)(1-T^2)) \\
+ QT(A+T)(1-T) + (A+T^2)(A + T)
\end{multline}
We have verified \eqref{eq: T33_shifting} for all $k \leq 90$.

\begin{rem}
There are distinct similarities between \eqref{eq: T2N_shifting} and \eqref{eq: T33_shifting}, both of which appear to be built from a few simple pieces. This hints at the existence of some nice closed form expression for $\mathcal{P}(n, dn)_{k, 1, \dots, 1}$. In some ways this is not so surprising; links of this form have braid presentations given by powers of full twist braids, and these braids are especially well-adapted to the usual methods of computing uncolored HOMFLY  homology. See e.g. \cite{Hog17b, EH19, GNR21}. We do not explore this line of inquiry any further here.
\end{rem}

\bibliographystyle{alpha}
\bibliography{Reduced_Homfly_bib.bib}

\end{document}